\documentclass[11pt,twoside]{amsart}
\usepackage{amsfonts}
\usepackage{amscd}
\usepackage{amssymb}
\usepackage{epsfig}
\usepackage{psfrag}

\title[A character formula for compact elements]{A character formula
for compact elements (the rank one case)}
\author{Jonathan Korman}
\address{University of Michigan, Ann Arbor, Michigan }
\address{Current Address: University of Toronto, Toronto, Ontario }
\email{jkorman@math.toronto.edu}

\begin{document}

\theoremstyle{plain}
\newtheorem{thm}{Theorem}
\newtheorem*{thm*}{Theorem}
\newtheorem{lem}[thm]{Lemma}
\newtheorem{cor}[thm]{Corollary}
\newtheorem{prop}[thm]{Proposition}
\newtheorem{claim}[thm]{Claim}
\newtheorem{eqn}[thm]{eqnarray}

\theoremstyle{definition}
\newtheorem{rem}[thm]{Remark}
\newtheorem{defn}[thm]{Definition}
\newtheorem*{defn*}{Definition}
\newtheorem{ex}[thm]{Example}

\maketitle

\begin{center} \bf Dissertation submitted to the University of
Michigan under the direction of Prof. Thomas C. Hales\\ 2002
\end{center}

\begin{abstract}
In their 1997 paper,~\cite{SS2}, Schneider and Stuhler gave a
formula relating the value of an admissible character of a
$p$-adic group at an elliptic element to the fixed point set of
this element on the Bruhat-Tits building. Here we give a similar
formula which works for compact elements. Elliptic elements have
finitely many fixed facets in the building but compact elements
can have infinitely many. In order to deal with the compact case
we truncate the building so that we only look at a bounded piece
of it. We show that for compact elements the (finite) information
contained in the truncated building is enough to recover all of
the information about the character. This works since the fixed
point set of a compact (non elliptic) element is periodic. The
techniques used here are more geometric in nature than the
algebraic ones used by Schneider and Stuhler. We recover part of
their result as a special case.
\end{abstract}

\setcounter{section}{-1}

\tableofcontents

\section{Introduction}

In this paper we use the geometry of the (semisimple) Bruhat-Tits
building to obtain a character formula for finitely generated
admissible representations of a connected reductive $p$-adic
group.

\subsection{}
Given an admissible representation $(\pi,V)$ of such a group, one
would like to understand its character $\Theta_\pi$, which is a
complex valued function on the group. Ideally one would like a
formula that expresses the value of the character at any given
element of the group.\\

Associated with a reductive $p$-adic group $G$ there is a
geometric space $X$, called the (semisimple) Bruhat-Tits building,
or affine building, of $G$. This building is
a polysimplicial complex endowed with a $G$-action.\\

As a geometric space, the building only encodes information about
its associated group, and not about a specific representation $V$
of that group. In order for it to contain information about such a
representation and its character, we have to consider a certain
sheaf on $X$ associated to $V$.
The action of the group $G$ on the building extends to this sheaf.\\

The formula for the value of characters (of finitely generated
admissible representations) on certain elements of $G$
is in terms of the action of these elements on $X$ and on this sheaf.\\

In 1997 Schneider and Stuhler ~\cite{SS2}, using algebraic
techniques, gave such a formula for elliptic elements -- elements
which have a finite number of fixed facets in the building. Using
different techniques, yet relying on some of their basic results,
we give a formula which works for compact elements -- which can
have an infinite number of fixed facets. In the elliptic
case, this gives a new proof of the Schneider-Stuhler result (for groups
of semisimple rank $1$).\\

The main idea is the use of certain truncation operators which
pick out finite subsets of the fixed-point set. These subsets
contain enough information from which to recover the character.
\subsection{}
Let $k$ be a $p$-adic field of characteristic zero. Let
$\mathbf{G}$ be a connected reductive algebraic group defined over
$k$ and $G=\mathbf{G}(k)$ its group of $k$-rational points. Denote
by $X=\mathcal{B}(\mathcal{D}\mathbf{G},k)$ the semisimple
Bruhat-Tits building; that is the Bruhat-Tits building associated
to the derived group, $\mathcal{D}G$, of $G$. Let $(\pi,V)$ be a
finitely generated admissible representation of $G$. Let
$U_{F}^{(e)}$ be the open-compact subgroups associated by
Schneider-Stuhler ~\cite{SS2} to the facets $F$ of $X$. These
subgroups satisfy the relations $U_{F'}^{(e)} \subset U_{F}^{(e)}$
whenever $F' \subset \overline{F}$. Consequently $V^{U_{F}^{(e)}}$
are finite dimensional subspaces of $V$ associated to the facets
$F$. The relations between the groups translate into relations
between
these vector spaces: $V^{U_{F}^{(e)}} \subset V^{U_{F'}^{(e)}}$.\\

In ~\cite{SS2}, Schneider and Stuhler prove:\\

$\bold{Theorem.}$ Let $G$ be a connected reductive group, $\gamma
\in G$ regular semisimple elliptic, and $(\pi,V)$ a finitely
generated admissible representation of $G$. Then there exists an
integer $e_0=e_0(V)$ which depends on $V$, such that for all $e
\geq e_0$ the character $\Theta_\pi$ can be expressed as:
\begin{eqnarray*}
\Theta_{\pi}(\gamma)= \sum_{q=0}^{d} \sum_{F(\gamma) \in
(X^{\gamma})_q} (-1)^q trace(\gamma,V^{U_{F}^{(e)}})
\end{eqnarray*}
Here $d$ is the dimension of $X$, $X^\gamma$ (a polysimplicial
complex, but not necessarily a subcomplex of $X$) is the fixed
point set of $\gamma$, $(X^{\gamma})_q$ is the set of $q$-facets
of $X^\gamma$ and $F(\gamma):= F \cap X^\gamma$ (a polysimplex).

$\bold{Remark.}$\label{stable/fixed} The summation in the formula
above should be understood to mean the sum over all of the
$\gamma$-stable facets in $X$.  When the action of $G$ on $X$
preserves the types\footnote{This happens for example if $G$ is
semisimple simply connected.} of the vertices in $X$, a facet is
$\gamma$-stable if an only if it is $\gamma$-fixed. Thus if this
is the case then the summation is over all the $\gamma$-fixed
facets and the above formula simplifies to:
\begin{eqnarray*}
\Theta_{\pi}(\gamma)= \sum_{q=0}^{d} \sum_{F \in X^{\gamma}_q}
(-1)^q trace(\gamma,V^{U_{F}^{(e)}})
\end{eqnarray*}

\subsection{}
Let $C_q:=C_c^{or}(X_{(q)};\gamma_e (V))$ be the vector space of
oriented $q$-chains with compact support (see ~\cite[II]{SS2}).
This is a \emph{smooth} representation of $G$. Write $T_g$ for the
action of $g \in G$ on the $C_q$'s. Consider the complex of
(smooth) $G$-modules:
\begin{gather}
 0 \longrightarrow C_d \overset{\partial}\longrightarrow \cdots
\overset{\partial}\longrightarrow C_0
\overset{\epsilon}\longrightarrow V \longrightarrow 0  \tag{$*$}
\end{gather}

For $V$ finitely generated admissible there exits an integer
$e_0=e_0(V)$ depending on $V$ such that for all $e\geq e_0$ the
complex $(*)$ is exact (~\cite{SS1} and ~\cite[Theorem II.3.1]{SS2}).\\

We would like to apply the \emph{Hopf trace formula} to the
operators $T_g$ acting on $(*)$. The $T_g$'s commute with the
boundary operators ($\partial$ and $\epsilon$), but they are not
of finite rank. In an attempt to address this issue we introduce
truncation operators $ Q_q^r$, ($r \geq 0$, a real number) on
$X_{(q)}$ which can be viewed as acting on $C_q$. Roughly, one can
think of these truncation operators as intersecting the building
with a ball of radius $r$ about some fixed point $o \in X^{\gamma}
$. An advantage of the $Q_q^r$'s is that they are of finite rank,
but unfortunately they do not commute with $\partial$. To fix
this, modify the truncation operators as follows.

Fix a vector space direct sum decomposition of $C_q$ (see section
~\ref{modtrunc} for details):
\begin{gather}
 C_q = B_q \oplus H'_q \oplus B'_{q-1} \tag{$\oplus$}
\end{gather}
Using such a decomposition it is possible
to define modified truncation operators, $\overline{Q}_q^r$ on
$C_q$. These modified operators have all of the desired
properties: they commute with $\partial$ (and $\epsilon$), they
have finite rank and they tend to the identity operator on $C_q$,
$\mathit{Id}_{C_q}$, as $r$
tends to $\infty$. Using these modified operators we obtain:\\

$\bold{Theorem.}$ Let $G$ be connected reductive, $V$ a finitely
generated admissible representation of $G$, $f \in C_c^\infty(G)$,
a locally constant, compactly supported function and $e\geq e_0$.
There exists a radius of truncation $r_0=r_0(f)$, which depends on
$f$, such that for all $r\geq r_0$ the trace of the operator
$\pi(f)$ can be expressed as:
$$trace\, \pi(f)= \sum_{q=0}^{d}  (-1)^q
trace(T_{f}\overline{Q}_q^r,C_q).$$

We would like to give a more geometric interpretation of this
formula in terms of the original truncation operators. With this
purpose in mind, assume that the direct sum decomposition is
$T_\gamma$-equivariant $(\gamma \in G^{cpt})$ and `nice'. Then
$trace(T_{\gamma}\overline{Q}_q^r,C_q)=trace(T_{\gamma}Q_q^r,C_q)$
and the above formula gives the main result (see Theorem
~\ref{main} for a more precise statement):
$$\Theta_\pi(\gamma)= \sum_{q=0}^{d} \sum_{F(\gamma) \in
(X^{\gamma} \cap X^r)_q} (-1)^q trace({\gamma},V^{U_{F}^{(e)}}).$$
Here $X^r$ is a finite subcomplex of $X$ called the truncated
building.\\

\subsection{}
For the following technical reasons we can prove this formula in
full generality only for groups of semisimple rank $1$. The direct
sum decomposition $(\oplus)$ is controlled by the following
truncated complex:
\begin{gather}
 0 \longrightarrow C_c^{or}(X_d^r ;V) \overset{\partial '}\longrightarrow
 \cdots
\overset{\partial '}\longrightarrow C_c^{or}(X_0^r ;V)
\overset{\epsilon '}\longrightarrow V   \tag{$\overline{*}$}
\end{gather}
 The exactness of this complex guarantees the existence of a
nice direct sum decomposition. Exactness of $(\overline{*})$ may
or may not depend on a parameter $e_r$. To prove the strong
version of the main result it is necessary for it not to depend on
$e_r$. We have the following:

\begin{itemize}
\item{ For $V=\mathbb{C}$ the trivial representation, the exactness
of $(\overline{*})$ follows from the contractibility of the
truncated building $X^r$, and is independent of $e_r$.}

\item{ For $G$ of semisimple rank $1$, we prove exactness (independently of $e_r$)
using the fact that $X$ has non-positive curvature, which implies
that the distance function on the building is strictly convex.}

\item{ For a general connected reductive group $G$, we prove
exactness of $(\overline{*})$ but the technique used is not
independent of $e_r$. Thus in this case we can only obtain a
weaker version of the main result.}
\end{itemize}

\subsection{}
The geometric techniques used here might suggest how one could
proceed with the case of non-compact elements. Such elements do
not have fixed points on $X$ but they do have fixed points on the
spherical building at infinity, $X^{\infty}$ , which is another
building associated with $G$. Finding a character formula for
non--compact elements is still an open problem.

\subsection*{Acknowledgments}
This paper is based on my thesis work under the supervision of Tom Hales at the University of Michigan.
None of this would have been possible without Tom's great optimism and endless
encouragements; it is a pleasure to thank him for all his help.

Special thanks to Ju-Lee Kim, to Jeff Adler and to my friends: Julia Gordon, 
Elliot Lawes and Joel Pitkin.

\newpage
\section{Notation}
\noindent $k$ a $p$-adic field of characteristic zero.\\
$k^\times$ the set of non zero elements of $k$.\\
$O$ the ring of integers of $k$.\\
$O^{\times}$ the set of units in $k$.\\
$\varpi$ a fixed generator for the maximal ideal in $O$.\\
$\overline{k}=O / \varpi O$ the residue field of $k$.\\
$\overline{k}^{\times}$ the set of non zero elements of $\overline{k}$.\\
$\omega :k^\times :\longrightarrow \mathbb{Z}$ the discrete
valuation normalized by $\omega(\varpi)=1$.\\
$G=\mathbf{G}(k)$ a connected reductive group.\\
$G^{reg}$ the set of regular semisimple elements in $G$.\\
$G^{ell}$ the set of regular semisimple elliptic elements in $G$.\\
$G^{cpt}$ the set of regular semisimple compact elements in $G$.\\
$\gamma$ a compact (sometimes also elliptic) element in $G$.\\
$X$ the semisimple Bruhat--Tits building of $G$.\\
$A$ a basic apartment in $X$.\\
$F$ a facet of $X$.\\
$\overline{F}$ the closure of the facet $F$ in $X$ (a
polysimplex).\\
$X^r$ a truncated building with truncation parameter $r$.\\
$X_q$ the $q$-dimensional facets of $X$.\\
$X_{(q)}$ the oriented $q$-dimensional facets of $X$.\\
$X^g$ the fixed point set of $g\in G$.\\
$(\pi,V)$ a finitely generated admissible representation.\\
$\Theta_\pi$ the character of $(\pi,V)$.\\
$\mathcal{O}_\gamma (f)$ the orbital integral of a function $f$
with respect to an element $\gamma$.\\
\newpage
\section{Preliminaries}
Let $k$ denote a $p$-adic field of characteristic zero, that is, a
finite extension of $\mathbb{Q}_p$ for some prime $p$. We will
denote by $O$ the ring of integers of $k$ and pick a generator
$\varpi$ for the maximal ideal in $O$. The residue field of $k$
will be denoted by $\overline{k}$.
Let $\mathbf{G}$ be a connected reductive group defined over $k$
and denote by $G$ the group $\mathbf{G}(k)$ of $k$-rational points
of $\mathbf{G}$, equipped with the natural locally compact
topology induced from that on $k$.\\
\subsection{The building}
A building is a polysimplicial complex which can be expressed as
the union of subcomplexes called apartments satisfying certain
axioms (see ~\cite[p.78]{Brown}). There are two types of
buildings: affine buildings and spherical buildings. The
apartments of affine buildings are Euclidean spaces and those of
spherical buildings are spheres. To a connected reductive group
$G$ one can associate at least three kinds buildings:
\begin{itemize}
  \item {} The semisimple Bruhat-Tits building of $G$.
   This is a building of affine type.
  \item {} The Bruhat-Tits building of $G$.
   This building is also of affine type.
  \item{} The spherical building (at infinity) of $G$.
  This is a building of spherical type.
\end{itemize}

\begin{rem}If $G$ is semisimple then the semisimple Bruhat-Tits building
and the Bruhat-Tits building are the same. For a reductive,
non-semisimple group $G$, the Bruhat-Tits building is a product of
the semisimple Bruhat-Tits building and an affine building
associated to the center of $G$.
\end{rem}

In this paper we will only use the semisimple Bruhat-Tits
building. We follow the review in ~\cite[I.1]{SS2} of the
construction of such a building.

\subsection{Review of the Semisimple Bruhat-Tits building}

We will use the following notation.\\

\noindent $G=\mathbf{G}(k)$ be a connected reductive group.\\
$S$ a maximal $k$-split torus in $G$.\\
$X^*(S):=\mathrm{Hom}_k(S,k^\times)$  the lattice of rational
characters of $S$.\\
$X_*(S):=\mathrm{Hom}_k(k^\times,S)$ the (dual) lattice of
rational co-characters of $S$.\\
$C:=Z(G)^\circ$ be the connected component of the center of $G$.\\
$X_*(C)$ the lattice of rational co-characters of $C$.\\
$Z$ the centralizer of $S$ in $G$.\\
$N$ the normalizer of $S$ in $G$.\\
$W:=N/Z$ the Weyl group.

\begin{defn*}
The underlying affine space of the real vector space
$$A:=(X_*(S)/X_*(C))\otimes \mathbb{R}$$
is called the \textit{basic apartment}.
\end{defn*}

The Weyl group $W$ acts by conjugation on $S$ which induces a
faithful linear action of $W$ on $A$.\\
Let $<\;,\;>\,:\,X_*(S)\times X^*(S) \longrightarrow \mathbb{Z}$
be the natural pairing; its $\mathbb{R}$-linear extension is also
denoted by $<\;,\;>$.\\
There is a unique homomorphism
$$\nu: Z \longrightarrow X_*(S)\otimes \mathbb{R}$$
characterized by
$$<\nu(g),\chi |_S>=-\omega(\chi(g))$$
for all $g\in Z$, and all characters $\chi$ of $Z$ (here $\omega$ is
the discrete valuation).\\
Using this homomorphism $g\in Z$ acts on $A$ by translations
$$gx:=x + \mbox{image of }\nu(g)\mbox{ in } A \quad x\in A.$$
This action of $Z$ on $A$ can be extended to an action of $N$ on
$A$. The $N$ action is compatible with the action of the Weyl
group $W$. Recall that there exists ~\cite[p.31]{Tits} a system of
affine roots $\Phi_{\mathit{af}}$ (which are certain affine
functions on $A$), and a mapping $\alpha \mapsto U_\alpha$ from
$\Phi_{\mathit{af}}$ onto a set of subgroups of $G$.

\begin{defn*}
Two points in $A$ are called equivalent if all affine roots have
the same sign on these two points; the corresponding equivalence
classes are called \textit{facets}. The facets of maximal
dimension $d$ are called \textit{chambers} (they are also the
connected components of the complements in $A$ of the union of
walls; a \textit{wall} is the zero set of an affine root). The
$0$-dimensional facets are the \textit{vertices}. The closure
$\overline{F}$ of a facet $F$ is a \textit{polysimplex} (in the
sense of algebraic topology). This gives the basic apartment a
polysimplicial structure.
\end{defn*}

Consider the following equivalence relation on $G \times A$:
$$(g,x) \sim (g',x') \mbox{ if there is an } n\in N \mbox{ such that }
  nx = x' \mbox{ and } g^{-1}g'n\in U_x.$$
We define $X:=G\times A/\sim$. It is easy to see that $G$ acts on
$X$. This action extends the action of $N$ on $A$ and the
polysimplicial structure of $A$ extends to $X$.

\begin{defn*}
The \textit{semisimple Bruhat-Tits building} of $G$ is the
polysimplicial $G$-complex $X$, also denoted by
$\mathcal{B}(\mathcal{D}\mathbf{G},k)$.
\end{defn*}

\subsection{Main properties of the building}
Let $X=\mathcal{B}(\mathcal{D}\mathbf{G},k)$ be the semisimple
Bruhat--Tits building of $G$; this is equivalent to saying that
$X$ is the Bruhat-Tits building of the derived group
$\mathcal{D}\mathbf{G}$ of $\mathbf{G}/k$.\\

We list the properties of the building which will be used in this
paper (see ~\cite{Brown}, ~\cite[1.1]{Moy} and ~\cite[I.1]{SS2}
for more details).

\begin{itemize}
  \item{} The building associated to $G$ is made up of
  \textit{apartments} which are glued together.
  Each apartment is a Euclidean space equipped with a polysimplicial
  structure and is isomorphic to $A$ (see ~\cite[p.10]{SS2}).

  \item {} The building is a $d$-dimensional locally finite
  polysimplicial complex, where $d=dim(A)$ is the semisimple
  $k$-rank of $G$. The group $G$ acts on $X$ polysimplicially.
  (If $\mathcal{D}\mathbf{G}$ is simple, as opposed to semisimple,
  then the building and the action are simplicial, as opposed to
  polysimplicial)

  \item{} Topologically, the building $X$ is a contractible space
  with a natural $G$-action.

  \item{} There is a natural metric $d(\cdot,\cdot)$ on $X$ with
  respect to which the action of $G$ is by isometries.

  \item {} Any two points (and even any two facets) in $X$ are
  contained in a common apartment.

  \item {} Any two points $x,y\in X$ are connected by a unique geodesic
  line segment, denoted $geod(x,y)$.

\end{itemize}
Here $X_q$ will denote the space of all $q$-facets of $X$ and
$X_{(q)}$ will denote the space of all oriented $q$-facets. Also
write $d=dim (X)$ for the dimension of $X$ as a locally finite
polysimplicial complex.

\begin{figure}[h]
\psfragscanon \centering  \psfrag{o}{$o$}
\epsfig{file=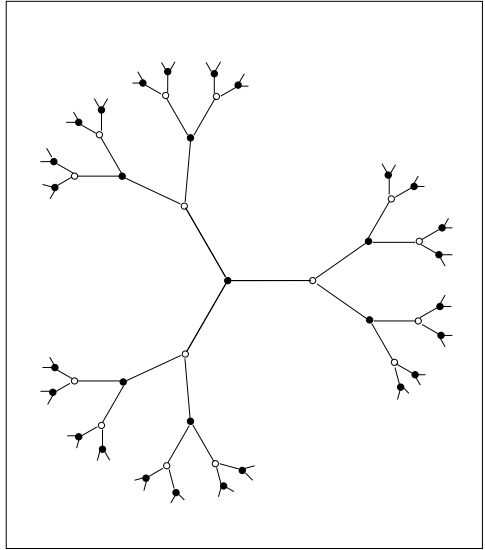,height=8cm}\caption{\label{tree0}The
(semisimple) Bruhat-Tits building of the group
$SL_2(\mathbb{Q}_2)$.}
\end{figure}

\begin{ex}
For $G=SL_2(\mathbb{Q}_{p})$, the semisimple Bruhat-Tits building
$X$ is a \emph{tree} with $p+1$ edges meeting at every vertex.
Since the group $SL_2$ is semisimple, its semisimple Bruhat-Tits
building and Bruhat-Tits building are the same. The apartments of
this building (tree) are one dimensional Euclidean spaces (lines).
The $0$-dimensional facets are the vertices and the
$1$-dimensional facets are the edges. See Figure ~\ref{tree0}.\\

Since $GL_2$ and $SL_2$ have the same derived group, their
semisimple Bruhat-Tits buildings are the same. For the actions of
$SL_2$ on this building there are two types of vertices; for the
action of $GL_2$ there is only one type\footnote{For explanation
of \textit{type} (in a special case), see ~\cite[p.30]{Brown}.}.
\end{ex}
\subsection{The character}\label{character}
We include Fiona Murnaghan's ~\cite{Murnaghan} explanation of the
character of an
admissible representation.\\

A representation $(\pi,V)$ of $G$ is called \textit{smooth} if
$$Stab_G(v):=\{g\in G \; |\; \pi(g)v=v\}$$ is open for every $v\in
V$.\\

Let $C_c^{\infty}(G)$ denote the space of complex valued, locally
constant, compactly supported function on $G$. Given $v\in V$ and
$f\in C_c^{\infty}(G)$, the function $g \mapsto f(g)\pi(g)v$
belongs to $C_c^{\infty}(G,V)$, compactly supported, locally
constant functions with values in $V$. Therefore, for each $f\in
C_c^{\infty}(G)$, we can define an operator $\pi(f)$ on $V$ as
follows:
$$ \pi(f)v = \int_G f(g)\pi(g)v\; dg, \qquad v\in V.$$
Here, $dg$ denotes a fixed Haar measure on $G$.\\

To define the character of $\pi$, we want to take the trace of
$\pi(f)$. Since $V$ is infinite dimensional in general this trace
is not defined. In order to make sense of the trace of $\pi(f)$ we
restrict the class of representations from smooth to admissible.

Recall that a smooth representation $(\pi,V)$ of $G$ is said to be
\textit{admissible} if for every open compact subgroup $K$ of $G$,
the space
$$V^K=\{v\in V \; |\; \pi(k)v=v \; \forall k\in K   \}$$
is finite dimensional. It can be shown that $\pi$ is admissible if
and only if $\pi(f)$ has finite rank for all $f$ in
$C_c^{\infty}(G)$. Thus for $\pi$ admissible we can talk about the
trace of the operator $\pi(f)$. Write
$$\Theta_\pi(f)=trace\, \pi (f), \qquad f\in C_c^{\infty}(G).$$
The distribution (linear functional on $C_c^{\infty}(G)$) defined
by $f \mapsto \Theta_\pi(f)$ is the character of $\pi$.

It is a theorem of Harish-Chandra (see ~\cite[p.2]{Murnaghan})
that for $\pi$ an admissible representation of finite
length\footnote{Recall that a finitely generated admissible
representation has finite length.}, there exits a locally
integrable function, also denoted by $\Theta_\pi$, on $G$, which
is locally constant on the set of regular elements in $G$, and
satisfies
$$\Theta_\pi(f)=\int_G f(g)\Theta_\pi(g)v\; dg, \qquad f\in
C_c^{\infty}(G).$$ That is, the distribution $\Theta_\pi$ is given
by integration against a function. The function $\Theta_\pi$ is
also called the character of $G$.\\

\subsection{Three types of elements}
Here $G^{reg}$ will denote the set of regular semisimple elements
in $G$. Recall that for an element $\gamma\in G^{reg}$, the
connected component of its centralizer,
$T:=C_G(\gamma)^{\circ}$, is a maximal torus in $G$.\\

We classify elements $\gamma \in G^{reg}$ into three
types\footnote{We could classify all elements in $G$, not just
those in $G^{reg}$, into these three types.} according to their
set of fixed points on the semisimple Bruhat--Tits building $X$.

\begin{defn} Let $\gamma \in G^{reg}$ be a regular semisimple
element.

We call $\gamma$ \textit{compact} if its fixed point set,
$X^\gamma$, is non-empty. Write $G^{cpt}$ for the set of regular
semisimple compact elements.\\
An equivalent characterization is:
\begin{itemize}
  \item {} $\gamma$ is contained in some subgroup of $G$ which is
  open, and compact modulo the center of $G$.
\end{itemize}

We call a compact element $\gamma$ \textit{elliptic}, if its fixed
point set, $X^{\gamma}$, is compact. Write $G^{ell}$ for the set
of regular semisimple elliptic elements.\\
Equivalent characterizations are:
\begin{itemize}
  \item {} $\gamma$ is not contained in any parabolic subgroup of $G$.
  \item {} $C_G(\gamma)^{\circ}$ is an elliptic torus (i.e.
         compact modulo the center of $G$).
\end{itemize}

We call $\gamma$ \textit{non-compact} if it has no fixed points on
the building $X$. (Sometimes non-compact
elements are also called \textit{hyperbolic}).\\
\end{defn}
\begin{rem}
Intuitively a non-compact element acts on $X$ by translations,
hence has no fixed points on $X$. Non-compact elements do have
fixed points on the spherical building at infinity $X^{\infty}$
(see ~\cite[VI.9]{Brown}). We will not deal with non-compact
elements here. Also note the following.

Let $G$ be a connected reductive group, and $(\pi,V)$ an
admissible representation of $G$. In ~\cite{Casselman}, Casselman
shows how to attach to a general element $g\in G^{reg}$, a
parabolic subgroup $P=MN$, such that $g\in M$ ($g$ is compact in
$M$). If $(\overline{\pi},V_N)$ denotes the representation of $M$
on the Jacquet module $V_N$, then it is shown in ~\cite{Casselman}
that $\Theta_\pi(g)=\Theta_{\overline{\pi}}(g)$. That is the
calculation of the value of a character at a general element $g$,
can be reduced to the case where $g$ is compact.
\end{rem}

\begin{figure}[h]
\psfragscanon \centering  \psfrag{A}{$A$}
\epsfig{file=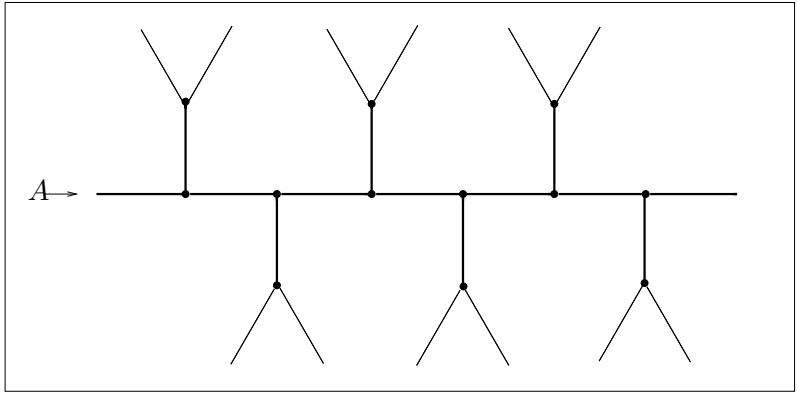,height=6cm}\caption{\label{tree1}The fixed
point set (in bold) of the compact non-elliptic element $\gamma$.}
\end{figure}

\begin{ex}\label{exnce}
The following elements in $\mathbf{GL_2}(k)$ provide examples of
the different types of regular semisimple elements. Here $\varpi$
is a uniformizer.
\begin{itemize}
  \item {elliptic:}  $\left( \begin{array}{cc}
               a   &   b\varpi \\
               b   &     a
           \end{array} \right)$, where $a \in k,\; b \in k^{\times}$ and
           $a^2 - b^2 \varpi \in k^{\times}$.\\
           The centralizer of this element is an elliptic torus.
  \item {compact, non elliptic:}   $\left( \begin{array}{cc}
               1   &   0 \\
               0   &   u
           \end{array} \right)$, where $u\in O^\times$ is a unit and $u \neq 1$.\\
           This element is contained in parabolic subgroups and it is easy
           to show (~\cite[Lemma 5.2]{Langlands}) that its fixed point
           set contains the basic apartment $A$.
  \item {non-compact:} $\left( \begin{array}{cc}
               1   &   0 \\
               0   &  \varpi^{n}
           \end{array} \right)$, $n \neq 0$ an integer.\\
           This element acts as a translation and so has no fixed
           points in $X$.
\end{itemize}

Since compact non-elliptic elements play a central role in this
paper, we take a closer look at such an element and its fixed
point
set (see ~\cite[chapet 5]{Langlands}).\\
Let $\gamma = \left(
\begin{array}{cc}
               1   &   0 \\
               0   &   u
           \end{array} \right)$
, where $u$ is a unit and $u \neq 1$.\\
Suppose $u$ has the form: $\;u= 1+\alpha_r\varpi^r+\cdots \in 1
+\varpi^r O^{\times}\subset O^{\times} \quad 0 \neq \alpha_r \in
\overline{k}.$ Then the points of $X$ fixed by $\gamma$ are
precisely those at a distance less than or equal to r from the
basic apartment $A$ ~\cite[Lemma 5.2]{Langlands}. For example, if
$r=1$, then the fixed point set looks like Figure ~\ref{tree1}.
\end{ex}

\subsection{\label{Ugroups}Some subgroups}
In ~\cite[I.1]{SS2}, Schneider and Stuhler attach to each facet
$F$ of the building the following subgroups:

\begin{itemize}
  \item {} The stabilizer of the facet:
  $P_F^{\dag}:=\{g\in G\; |\;gF=F\}$.
  \item {} The `fixer' of the facet: $P_F:=\{g\in G \;|\; gz=z \;
  \mbox{for all} \; z\in F\}$.\\
  $P_F^{\dag}$ and $P_F$ are open subgroups of $G$;
  $P_F$ is a subgroup of finite index in $P_F^{\dag}$.
  \item {} A filtration of subgroups: $U_F^{(e)}$,
  parameterized by integers $e \geq 0$. These are open compact in $G$
  and normal in $P_F^{\dag}$.
\end{itemize}

Sometimes it will be convenient to attach such groups to any point
in the building. Let $x\in X$ be any point in $X$ and let $F$ be
the unique facet containing $x$. Define $P_x^{\dag}=P_x=\{g\in G\;
|\;gx=x\}$ and $U_x^{(e)}=U_F^{(e)}$.

We summarize some facts about the groups $U_F^{(e)}$ which will be
used later. See ~\cite[Chapter I]{SS2} and ~\cite[Lemma
1.28]{Vig}. Here $e \geq 0$ is an integer and $F$ is any facet of
the building $X$.
\begin{description}

\item[(U1)] $U_F^{(e)}$ are open compact subgroups of $G$.
See ~\cite[p.13]{SS2}.\\

\item[(U2)] $U_F^{(e)} \lhd P_F^{\dagger}$. See ~\cite[p.21]{SS2}.\\

\item[(U3)] $U_F^{(e)}$, for $e \geq 0$, form a fundamental system of neighborhoods of the identity element, $id$,
in $G$. See ~\cite[I.2 Corollary 9]{SS2}.\\

\item[(U4)] $U_{F'}^{(e)} \subset U_F^{(e)}$ for any two facets
$F',F$ in $X$ such that $F' \subset \overline{F}$.
See ~\cite[I.2 Proposition 11]{SS2}.\\

\item[(U5)] For any two vertices $x,y \in \overline{F}$ and for any $e,e' \geq 0$,
the subgroups $U_x^{(e)}$ normalize the subgroups $U_y^{(e')}$.
See ~\cite[end of section I.2]{SS2}.\\

\item[(U6)] $$U_F^{(e)}= \prod_{x \; vertex \; in \; \overline{F}}
U_x^{(e)}$$ for any facet $F$ in $X$ and any ordering of the
factors on the right hand side. See ~\cite[I.2 Proposition 11]{SS2}.\\

\item[(U7)] Fix two different points $x$ and $y$ in $X$. Recall
that $geod(x,y)$ denotes the (closed) geodesic joining $x$ and
$y$. In an apartment $A$ containing both $x$ and $y$ this geodesic
can be realized as $$geod(x,y)=\{ (1-t)x + ty \; | \; 0 \leq t
\leq 1 \}.$$ For any point $z \in geod(x,y)$ and for any $e  \geq
0$ we have:
$$U_z^{(e)} \subset U_x^{(e)}U_y^{(e)}.$$ See ~\cite[I.3
Proposition 1]{SS2} for $x$ a special vertex; remark after proof
of this proposition and ~\cite[Lemma 1.28]{Vig} for $x$ any
vertex.
\end{description}

\subsection{Representations as coefficient systems}
Let $\mathrm{Alg}(G)$ be the category of smooth representations of
$G$. In ~\cite[II.2]{SS2}, Schneider and Stuhler define the
following
objects.\\

\begin{defn*} A \textit{coefficient system} (of complex vector
spaces) $\underline{\underline{V}}$ on the Bruhat-Tits building
$X$ consists of:
\begin{itemize}
  \item {} complex vector spaces $V_F$ for each facet $F \subset
  X$, and
  \item {} linear maps $r_{F'}^F: V_F \rightarrow V_{F'}$ for each
  pair of facets $F' \subset \overline{F}$ such that $r_{F}^F=id$
  and $r_{F''}^F= r_{F''}^{F'} \circ r_{F'}^{F}$ whenever
  $F'' \subset \overline{F'}$ and $F' \subset \overline{F}$.
\end{itemize}
\end{defn*}
Coefficient systems form a category denoted by
$\textrm{Coeff}(X)$.\\

Fix an integer $e \geq 0$. For any representation $V$ in
$\textrm{Alg}(G)$ we have the coefficient system
$\underline{\underline{V}}:={(V^{U_F^{(e)}})}_F$. Write $\gamma_e:
\textrm{Alg}(G) \longrightarrow \textrm{Coeff}(X)$ for the
functor: $V \longrightarrow {(V^{U_F^{(e)}})}_F$.\\

Recall that $X_{(q)}$ denotes the space of all oriented $q$-facets
of $X$. We will denote an oriented facet by $(F,c)$, where $c$ is
an orientation of $F$ (see ~\cite[II.1]{SS2}). For any $F\in
X_{q+1}$ and any $F' \subset \overline{F}$ the map
$\partial_{F'}^F$ (~\cite[pp. 28--29]{SS2}), takes an orientation
$c$ on $F$ and returns (the induced) orientation
$\partial_{F'}^F(c)$, on $F'$. Often we will abuse the notation
and write $F$, instead of $(F,c)$, for an oriented facet.

\begin{defn*} For any $d \geq q \geq 0$ the space of
\textit{oriented q-chains of compact support with values in the
coefficient system} $\gamma_e(V)$ is:\\
$C_q:=C_c^{or}(X_{(q)} ;\gamma_e(V)):= \mathbb{C}-$vector space of all
maps $\omega:X_{(q)} \rightarrow V$ such that
\begin{itemize}
  \item {} $\omega$ has finite support,
  \item {} $\omega((F,c))\in V^{U_F^{(e)}}$, and
  \item {} if $q \geq 1$, $\omega((F,-c))=-\omega((F,c))$ for all
  $(F,c)\in X_{(q)}$.
\end{itemize}
\end{defn*}

The group $G$ acts \textit{smoothly} on these spaces via
$$(g\omega)((F,c)):=g(\omega((g^{-1}F,g^{-1}c))).$$
There is a natural boundary map
\begin{eqnarray*}
\partial:C_c^{or}(X_{(q+1)}
;\gamma_e(V)) & \longrightarrow & C_c^{or}(X_{(q)} ;\gamma_e(V))\\
\omega & \mapsto& ((F',c')\rightarrow \underset{(F,c) \in
{X_{(q+1)}} \atop{F' \subseteq \overline{F} \atop
{\partial_{F'}^F(c)=c' }}}\sum \omega ((F,c))).
\end{eqnarray*}

which satisfies: $\partial \circ \partial = 0$.

We obtain the following augmented $G$-equivariant chain complex:
\begin{gather}
 0 \longrightarrow C_d \overset{\partial_d}\longrightarrow \cdots
\overset{\partial_1}\longrightarrow C_0
\overset{\epsilon}\longrightarrow V \longrightarrow 0 \tag{$*$}
\end{gather}
were the augmentation map is given by
\begin{eqnarray*}
\epsilon: C_c^{or}(X_{(0)};\gamma_e(V)) & \longrightarrow & V\\
\omega & \mapsto & \underset{F \in {X_{(0)}}}\sum
\omega(F).\\
\end{eqnarray*}

Let $V\in\textrm{Alg}(G)$. In ~\cite[II.2-3]{SS2}, Schneider and
Stuhler prove the following:

\begin{itemize}
  \item {} If $V$ is admissible then
  $C_c^{or}(X_{(q)};\gamma_e(V))$ are finitely generated smooth
  $G$-modules.
  \item {} If $V$ is finitely generated then there exists
  $e_0=e_0(V)$ such that the complex $(*)$ is exact for all
  $e \geq e_0$.
\end{itemize}


For the rest of this paper $V$ will be a finitely generated
admissible representation of $G$, and $e_0$ large enough so that
the complex $(*)$ is exact.

\section{Overview of the Schneider-Stuhler result}
The main result of ~\cite{SS2} in which we are interested, can be
formulated as follows\footnote{This result is ~\cite[III.4 Lemma
10]{SS2} combined with \cite[III.4 Proposition 16]{SS2}. A
cohomological interpretations of this result is given by the
Hopf-Lefschetz type trace formula of \cite[IV.1 Proposition
5]{SS2}.}:

For $G$ a connected reductive group, $(\pi,V)$ a finitely
generated admissible representation of $G$, $\gamma\in G^{ell}$ a
regular semisimple elliptic element, and $e \geq e_0(V)$, we have:
\begin{eqnarray}\label{chi}
\Theta_\pi(\gamma) = \sum_{q=0}^{d}  \underset{F(\gamma) \in
{(X^\gamma)}_q} \sum (-1)^q trace(\gamma; V^{U_F^{(e)}})
\end{eqnarray}

Here $\Theta_\pi$ is the character (function) of $(\pi,V)$
evaluated on the regular element $\gamma$, $F(\gamma):=F \cap
X^\gamma$, and ${(X^\gamma)}_q$ are the $q$-facets of $X^\gamma$.
As stated in ~\cite[p.635]{Kottwitz}, $F(\gamma)$ is a polysimplex
and $X^\gamma$ is a polysimplicial complex.

\begin{rem}Since $\gamma F=F$ for a facet $F$ if and only if $F
\cap X^\gamma \neq \emptyset$, the sum over $F(\gamma) \in
{(X^\gamma)_q}$, over all $q$, is the same as the sum over all the
$\gamma$-stable facets $F \in X$.
\end{rem}

We now loosely sketch (our interpretation) of the proof given by
Schneider and Stuhler in ~\cite{SS2} of this main result. We will
sometimes use the phrase `for sufficiently large $e$' to mean `for all $e \geq e_0$'.\\

There are two main steps:\\
\begin{description}

\item[Step 1]\label{step1} Establish the form of the character formula
in terms of the fixed point set $X^\gamma$ for all \textit{large}
$e$, i.e. $e \geq e(f)$, where $e(f)$ is a constant which depends
on $f$:
$$\int_G f(g) \Theta_\pi(g)\, dg = \int_G f(g)
\left\{ \sum_{q=0}^{d} \underset{F(g) \in {(X^g)}_q} \sum (-1)^q
trace(g; V^{U_F^{(e)}}) \right\}\, dg $$ Here $f\in
C^{\infty}_c(G)$ has support in $G^{ell}$ so that $X^g$ is compact
and the sum on the right hand side is finite.\\

\item[Step 2] Show that for each fixed $g\in G^{ell}$ the following expression is independent of $e$
(for $e \geq e_0$):
$$\sum_{q=0}^{d} \underset{F(g) \in {(X^g)}_q} \sum (-1)^q trace(g;
V^{U_F^{(e)}})$$

\end{description}
Once independence of $e$ has been shown in {\bf Step 2}, we go
back to the formula of {\bf Step 1}, whose validity now holds for
all $e \geq e_0$. We choose $f$ to have support on a small
neighborhood of $\gamma$, small enough so that both the character
and the alternating sum above are constant on it. Thus we can
eliminate $f$ from the equation and obtain formula ~(\ref{chi})
to hold at an element $\gamma\in G^{ell}$.\\

The proof of {\bf Step 1} is essentially included in Lemmas 13 and
14 of ~\cite[III.4]{SS2}. We now give a short derivation of the
same result using elementary facts.

For simplicity of exposition we assume now that all stable facets
are actually fixed facets (see Remark ~\ref{stable/fixed} and the
discussion after Example~\ref{PGL_2}).\\

Fix $\gamma \in G^{cpt}$ and let $F\in X^\gamma$ be a
$\gamma$-fixed facet. Since $\gamma \in P_F^{\dagger}$ it acts on
$V^{U_F^{(e)}}$. Recall Harish-Chandra's formula (section
~\ref{character}):
$$\Theta_\pi(f)= \underset{G}\int f(g) \Theta_\pi(g)\, dg  \qquad
f\in C^{\infty}_c(G)$$

Let $K:=U_F^{(e)}$, an open compact neighborhood of the identity,
and choose $e(\gamma)$ large enough so that for all $e \geq
e(\gamma)$, (the locally constant function) $\Theta_\pi$ is
constant on $\gamma K$: $\Theta_\pi (\gamma K) = \Theta_\pi
(\gamma)$. Use $f=\frac {1_{\gamma K}}{vol(K)}$ in the
Harish-Chandra's formula to obtain:

\begin{eqnarray*}
\Theta_\pi (\gamma) & = & \int_K \frac {1}{vol(K)} \Theta_\pi (\gamma k)
\,dk\\
 & = & trace(\pi (\frac {1_{\gamma K}}{vol(K)}); V)\\
 & = & trace(\pi (\gamma) \pi (K) ;V)\\
 & = & trace(\pi (\gamma) ;V^K)\\
 & = & trace(\pi (\gamma) ;V^{U_F^{(e)}})
\end{eqnarray*}

\begin{rem} This says that for $\gamma \in G^{cpt}$ (so in particular
for $\gamma \in G^{ell}$), all the information about the value of
the character at $\gamma$ is contained in (the vector space
$V^{U_F^{(e)}}$ above $F$ for) any fixed facet $F$ of $\gamma$.
But to extract this information we are forced to increase the
parameter $e$. Unfortunately we have no clear control over the way
in which $e$ increases.
\end{rem}

We now recover the formula of {\bf Step 1}. Recall that for any
$\gamma\in G^{cpt}$ its fixed point set, $X^\gamma$, is
contractible: for two $\gamma$-fixed points $x,y$ their geodesic
$geod(x,y)$ is also $\gamma$-fixed, hence $X^\gamma$ is
(geodesically) contractible; and that for $\gamma\in G^{ell}$ this
fixed point set is compact (in fact it is a finite polysimplicial
complex).

Let $\gamma\in G^{ell}$. Choose $e=e(\gamma)$ large enough so that
$\Theta_\pi (\gamma U_F^{(e)}) = \Theta_\pi (\gamma)$, for all
(finitely many) facets $F \in X^\gamma$. Thus
$\Theta_\pi(\gamma)=trace(\gamma; V^{U_F^{(e)}})$, for all facets
$F\in X^\gamma$. Using the contractibility of $X^\gamma$, the
alternating sum below collapses and we recover formula
~(\ref{chi}), which is essentially the formula of {\bf Step 1}
(here $Euler(X^\gamma)$ is the Euler characteristic of
$X^\gamma$):
\begin{eqnarray*}
\sum_{q=0}^{d}  \underset{F\in X_q^\gamma}\sum (-1)^q
trace(\gamma; V^{U_F^{(e)}}) &=& Euler(X^\gamma)
trace(\gamma; V^{U_F^{(e)}})\\
&=& trace(\gamma; V^{U_F^{(e)}}) = \Theta_\pi (\gamma)
\end{eqnarray*}

\begin{rem}
The only facts that went into the derivation of this formula are
the local constancy of the character $\Theta_\pi$ and the
contractibility of fixed point sets. In this sense our derivation
is elementary. This means that the form of the character formula
is not too surprising and that it is in showing that this formula
is independent of the parameter $e$ ({\bf Step 2}), where most of
the effort is expanded.
\end{rem}

We now sketch the main ideas in {\bf Step 2}.\\
Recall the exactness of the resolution $(*)$ for sufficiently
large $e$ (i.e. $e \geq e_0$), and the fact that the $G$-modules
$C_q=C_c^{or}(X_q ;\gamma_e(V))$
are finitely generated. Both facts are proved in ~\cite{SS2}. \\

Consider the expression:
$$\sum_{q=0}^{d} (-1)^q dim \; H^q (\textrm{Hom}_G(C_q;V'))$$

Here $V'=(\pi',V')$ is any irreducible admissible representation
of $G$. Since $C_q$ are finitely generated as $G$-modules, the
vector spaces $\textrm{Hom}_G(C_q;V')$ are finite dimensional and
so the expression above makes sense. Since the resolution is exact
for all $e\geq e_0$, this expression is independent of $e$ (for
$e\geq e_0$).

Using the {\it Hopf trace formula} and a formal fact about
representations of finite groups, the expression above is related
to the \textit{Euler-Poincare function}\footnote{For definition of
the Euler-Poincare function see ~\cite[p.45]{SS2}} $f^{V,e}_{EP}$:
\begin{eqnarray*}
\sum_{q=0}^{d} (-1)^q dim \; H^q (\textrm{Hom}_G(C_q;V')) &=&
\sum_{q=0}^{d} (-1)^q dim (\textrm{Hom}_G(C_q;V'))\\ &=&
trace(\pi'(f^{V,e}_{EP});V')
\end{eqnarray*}

The left hand side being independent of $e$ implies that the right
hand side is also independent of $e$ for all irreducible
admissible representations $V'$.\\
Recall {\it Kazhdan's Density Theorem}~\cite[p.29]{Kazhdan}:

\begin{thm*}
Let $G$ be a connected reductive $p$-adic group. Let $f\in
C_c^{\infty}(G)$. Suppose that $trace(\pi'(f);V')=0$ for all
irreducible admissible representations $(\pi',V')$ of $G$. Then
$\mathcal{O}_\gamma (f)=0$ for all strongly regular\footnote{A
semisimple element is regular if its centralizer has the lowest
possible dimension. A semisimple regular element is strongly
regular if its centralizer is connected. In general an element is
regular/stronly regular if its semisimple part is regular/strongly
regular.} semisimple elements $\gamma \in G$.
\end{thm*}
Since $trace(\pi'(f^{V,e}_{EP});V')$ is independent of $e$ we
obtain for any other $e'$ :
$$0=trace(\pi'(f^{V,e}_{EP});V')-trace(\pi'(f^{V,e'}_{EP});V')
=trace(\pi'(f^{V,e}_{EP}-f^{V,e'}_{EP});V')$$ So by applying
Kazhdan's Density Theorem, we obtain:
$$0=\mathcal{O}_\gamma (f^{V,e}_{EP}-f^{V,e'}_{EP})=
\mathcal{O}_\gamma (f^{V,e}_{EP}) -  \mathcal{O}_\gamma
(f^{V,e'}_{EP}).$$

That is, the orbital integral $\mathcal{O}_\gamma (f^{V,e}_{EP})$
is also independent of $e$.\\

A formal calculation, ~\cite[III.4 Lemma 10]{SS2}, equates this
orbital integral and the alternating sum in equation ~(\ref{chi}):
$$\mathcal{O}_\gamma(f^{V,e}_{EP})=\sum_{q=0}^{d}  \underset{F\in X_q^\gamma}\sum (-1)^q
trace(\gamma; V^{U_F^{(e)}})$$ showing that the alternating sum is
independent of $e$.

\section{Truncated buildings and truncation operators}
Fix a point $o \in X$ and call it the \emph{origin}.\\
In this section we define a family of subsets, $X^r$ ($r \geq 0$,
a real number), of the building $X$. Each of these subsets will be
referred to as a \emph{truncated building} (with center $o$).
Using such a truncated building $X^r$, we will define
\emph{truncation operators} $Q^r_q$ on the $C_q$'s.\\

\textbf{Warning}: Do not confuse the truncated building $X^r$ with
fixed point sets (such as $X^g$ or $X^\gamma)$ or with $X_q$ for
that matter.\\

Let $B(o,r):=\{x\in X \;|\; d(o,x) \leq r\}$ be the closed ball of
radius $r$ about the point $o$. Here we denote by $d(\cdot ,\cdot
)$ the distance function on the building.
\begin{defn}
Let $S$ be any subset of the building. We define two operations on
$S$:
\begin{enumerate}
\item{ $cnvx(S):= \emph{convex hull of } S$.\\
We call a subset of the building \emph{convex} if for any two
points $x,y$ in the subset, the geodesic $geod(x,y)$ is also in
the subset. The convex hull of $S$ is the (unique) smallest convex
subset of the building containing $S$.\\
Note that if $o \in S \cap X^{\gamma}$ and $S$ is
$\gamma$-invariant, then $o \in cnvx(S)$ and $cnvx(S)$ is also
$\gamma$-invariant.}\\

\item{ $simp(S):= \emph{smallest subcomplex of the building containing S}$.\\
Note that, since the action of $\gamma$ is simplicial, if $o \in S
\cap X^\gamma$ and $S$ is $\gamma$-invariant, then $o \in simp(S)$
and $simp(S)$ is also $\gamma$-invariant.}
\end{enumerate}
\end{defn}
Now start with the subset $S=B(o,r)$ and apply the operations (1)
and (2) above consecutively to obtain the sequence of increasing
subsets:
\begin{gather}
 S \subseteq cnvx(S) \subseteq simp(cnvx(S)) \subseteq
\cdots \tag{$\dagger$}
\end{gather}

\begin{claim} This sequence stabilizes after a finite number of
terms.
\end{claim}
\begin{proof}
We assume the results of section ~\ref{subcomplex} where it is
shown (see Proposition ~\ref{ecnvx}) that there exist finite
polysimplicial subcomplexes of the building which are convex and
are arbitrary large (in the sense that they contain any given
ball). Starting with $B(o,r)$, fix such a subset $E$,
($E=X^{\gamma_{m(r)}}$ in the notation of ~\ref{subcomplex}),
containing it. Note that if $S \subseteq E$ then $cnvx(S)
\subseteq E$ and also $simp(S)\subseteq E$. Thus the terms in the
sequence $(\dagger)$ can never leave the finite subcomplex $E$ and
hence the sequence stabilizes.
\end{proof}
We denote by $X^r$ the stable terms in the sequence $(\dagger)$ and call
it the \emph{truncated building (with parameter r)}.\\

The main properties of the truncated building $X^r$ are:
\begin{itemize}
\item{$o \in B(o,r) \subseteq X^r$}.
\item{$X^r$ is convex}.
\item{$X^r$ is a finite subcomplex. (In particular it is compact)}
\end{itemize}

\begin{rem} If $B(o,r)$ contains a chamber, (e.g. $r$ not too
small), then $X^r$ is a union of maximal dimensional
polysimplices. (See Lemma ~\ref{fatX^g}).
\end{rem}

\begin{ex}\label{SL2ex} For $G=SL_2(\mathbb{Q}_{p})$, $X$ is a \emph{tree}
with $p+1$ edges meeting at every vertex. In this case the
truncated building $X^r$, for $r$ an integer, \textit{is} the
closed ball $B(o,r)$. See Figure ~\ref{trunctree}.\\
\begin{figure}[h]
\psfragscanon \centering  \psfrag{o}{$o$}
\epsfig{file=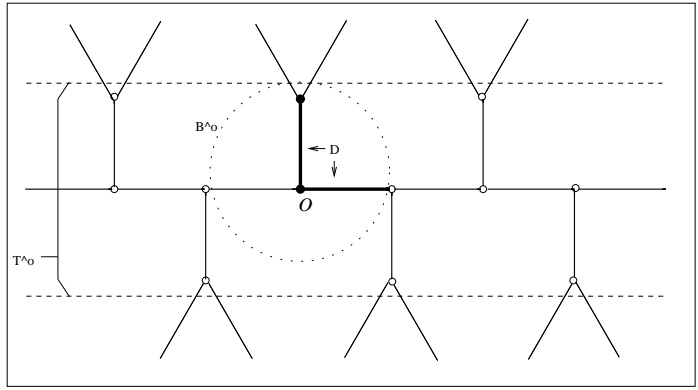,height=6cm}\caption{\label{trunctree}The
truncated building $X^r$ (part of building inside the dotted
circle) for $SL_2$ with $r=1$ and $r=2$.}
\end{figure}
To see this recall (see ~\cite[IV.3]{Brown}) that the distance
function, $d(o,\cdot)$ on $X$, is a geodesically convex function
\footnote{A function $f$ on $X$ is geodesically convex if for any
geodesic $geod(x(0),x(1)):=\{x(t)=tx+(1-t)y\;|\; 0 \leq t \leq 1
\}:\; f(x(t)) \leq tf(x(0))+(1-t)f(x(1)) \mbox{ for all } 0 \leq t
\leq 1$.}.

Hence the ball $B(o,r)$ is a convex set, and so applying $cnvx$ to
it does nothing. Now for $G=SL_2$, the ball $B(o,r)$ is already
simplicial, so applying $simp$ to it also does nothing.
\end{ex}

The truncated building $X^r$ induces a decomposition of each
oriented $q$-skeletons as a disjoint union: $$X_{(q)} =
X_{(q)}^{in(r)} \amalg X_{(q)}^{out(r)},$$ where
$X_{(q)}^{in(r)}:=\{F\in X_{(q)} \;|\; F \subset X^r
\}$ and $X_{(q)}^{out(r)}:=\{F\in X_{(q)} \;|\; F \not\subset X^r\}$.\\
Let $C_q^{in(r)}:=\mbox{the oriented q-chains in } \, C_q
\mbox{ supported on } X_{(q)}^{in(r)}$ and\\
let $C_q^{out(r)}:=\mbox{the oriented q-chains in } \, C_q
\mbox{ supported on } X_{(q)}^{out(r)}$.\\
Note that $C_{q}^{in(r)} = C_q(X^r)$.\\

We obtain a vector space direct sum decomposition:
$$C_q=C_q^{in(r)} \oplus C_q^{out(r)}$$
Define the \emph{truncation operator} $Q_q^r$ on $C_q$ to be the
projection of $C_q$ onto $C_q^{in(r)}$. Thus the truncation
operator $Q_q^r$ takes a map $\omega\in
C_q^{or}(X_{(q)};\gamma_e(V))$ and returns a map $Q_q^r\omega\in
C_q^{or}(X_{(q)};\gamma_e(V))$ supported on $X^r$.\\

Note the following properties of the truncation operators $Q_q^r$:
\begin{itemize}
\item{$Q_q^r$ is a finite rank operator.\\ This follows from the
fact that its range $C_q^{in(r)}$ is a finite dimensional vector
space}.
\item{The $Q_q^r$'s do not necessarily commute with the boundary
maps $\partial$.}
\end{itemize}

{\bf Notation:} Write $Q^r$ for the sequence of operators $(Q^r_d,
\cdots ,Q^r_0)$, where each operator $Q^r_q$ acts on the
corresponding vector space $C_q$.\\

We would like to modify the truncation operators so that they
\emph{do} commute with the boundary maps $\partial$. This
is done in section ~\ref{modtrunc}.\\

\subsection{\label{subcomplex}Existence of certain finite subcomplexes}
We now prove the existence of finite subcomplexes of the building
which are convex and are arbitrary large (in the sense that they
contain any given ball).\\
This section is independent of previous results.
\begin{lem}\label{fatX^g} Let $\gamma \in G=\mathbf{G}(k)$ be a compact
element of $G$. Then its fixed point set, $X^\gamma$, is a convex
subset of the building $X$. If $\gamma$ is elliptic, then
$X^\gamma$ is also compact.
\end{lem}

\begin{proof}
Convexity follows from the fact that for any two points $x,y \in
X$ fixed by $\gamma$, their geodesic $geod(x,y)$ is also fixed by
$\gamma$. Compactness is proved in \textit{Lemma 1} of
~\cite{Rog}.
\end{proof}

\begin{lem} For any compact element $\gamma \in G$, if its fixed point
set $X^\gamma$ contains a chamber $C$, then $X^\gamma$ is the
union of polysimplices of maximal dimension; equivalently
$X^\gamma$ is the union of the closure of its chambers.
\end{lem}
\begin{proof}
Pick any element $x \in X^\gamma$. Choose an apartment $A$
containing both $x$ and the chamber $C$. Recall that any apartment
is a Euclidean space (over $\mathbb{R}$) with the usual topology.
Let $cone(C,x)$ be the \emph{cone on $C$ with vertex $x$}. Since
both $C$ and $\{x\}$ are subsets of $X^\gamma$, so is $cone(C,x)$.
Let $star^\circ (x)$ denote the union of all the chambers $C'$ in
$A$ which contain $x$ in their closure. Both $cone(C,x)$ and
$star^\circ (x)$ are open subsets of the Euclidean space $A$ and
so is their non-empty intersection. Let $C' \subseteq star^\circ
(x)$ be a chamber which intersects $cone(C,x)$. This intersection
is fixed under $\gamma$ and it being an open subset of $C'$,
forces $C'$ to also be fixed by $\gamma$. Thus $x \in
\overline{C'}$.
\end{proof}

\begin{rem} For an arbitrary elliptic element
$\gamma$, its fixed point set $X^\gamma$ need not be a subcomplex
of the building.
\end{rem}
\begin{ex}\label{PGL_2} In $G=PGL_2(\mathbb{Q}_{2})$, which has the same
Bruhat-Tits building as $SL_2(\mathbb{Q}_{2})$, the fixed point
set of the elliptic element
\[ \gamma =\left[ \left( \begin{array}{cc}
      0      & 1 \\
      \varpi & 0
          \end{array} \right) \right]            \]
is a single point: the center of a standard chamber. This is not a
subcomplex of the building.
\end{ex}
For this reason we will restrict ourselves\footnote{Recall that
the groups $G$, $\widetilde{G}$, and $\mathcal{D}G$ all have the
same semisimple Bruhat-Tits building and their actions are
compatible.} to elliptic elements $\gamma \in
\mathcal{D}G=\mathcal{D}\mathbf{G}(k)$ which are coming from the
simply connected cover\footnote{By definition the simply connected
cover $\widetilde{G}$ of a reductive group $G$ is the simply
connected cover of its derived group $\mathcal{D}G$.}
$\widetilde{G}=\mathbf{\widetilde{G}}(k)$ of $G$. That is, let
$\widetilde{\gamma}$ be an elliptic element in $\widetilde{G}$ and
consider its image $\gamma$ under the natural map from
$\widetilde{G}$ to $\mathcal{D}G$. (Recall that the projection map
$\pi:\widetilde{\mathbf{G}}\longrightarrow \mathcal{D}\mathbf{G}$,
which is a map on the level of algebraic groups, induces a map
$\pi:\widetilde{\mathbf{G}}(k) \longrightarrow
\mathcal{D}\mathbf{G}(k)$ on the level of $k$-rational points and
that the image $\pi ( \widetilde{\mathbf{G}}(k))$ is a subgroup of
$\mathcal{D}\mathbf{G}(k)$ of finite index.) The action of the
simply connected cover $\widetilde{G}$ preserves the type (See
~\cite[2.5]{Tits}) of each vertex (do not confuse type of a vertex
with type of an element). Since any facet contains at most one
vertex of each type, any element $\gamma\in \widetilde{G}$ which
stabilizes a facet must fix this facet. That is, if an element
$\gamma\in \widetilde{G}$ fixes a point $x \in X$, then it fixes
every point of the unique facet in which $x$ lies.  Thus the fixed
point sets
$X^\gamma$ of such elements are subcomplexes.\\

Let $\tilde{\gamma}_m \in \widetilde{G}$ be a sequence of regular
elliptic elements going to $id \in \widetilde{G}$. Denote their
images in $\mathcal{D}G$ by $\gamma_m$.
\begin{lem} For each $r \geq 0$, there is a
positive integer $m=m(r)$ such that $\gamma_m$ fixes all points of
$B(o,r)$.
\end{lem}
\begin{proof}
Consider the open subgroup $U:=\{g\in \mathcal{D}G\; |\; gx=x \;
\mbox{ for all } x\in B(o,r) \}$. For large enough $m=m(r)$, the
element $\gamma_m$ is inside $U$ and so fixes all the points of
$B(o,r)$
\end{proof}

Combining the last two lemmas we obtain:
\begin{prop}\label{ecnvx} For any ball $B(o,r)$ the subset $X^{\gamma_{m(r)}}$
contains the ball. It is a convex, finite subcomplex of $X$. If
$B(o,r)$ contains a chamber, $X^{\gamma_{m(r)}}$ is a union of
maximal dimensional polysimplices.
\end{prop}

\section{\label{modtrunc}Modified truncation operators}
Some of the ideas in this section were inspired by ~\cite{AB1} and
~\cite{AB2}.\\

Fix direct sum decompositions of the vector spaces $C_q$ ($d \geq
q \geq 0$):
\begin{gather}
C_q = Z_q \oplus B'_{q-1} = B_q \oplus H'_q \oplus B'_{q-1}
\tag{$\oplus$}
\end{gather}

Here $Z_q$ are the $q$-cycles and $B'_{q-1}$ is a complement to
$Z_q$ inside $C_q$. The spaces $B_q:=\partial (C_{q+1})$ are the
$q$-boundaries, and so $B'_{q-1} \cong B_{q-1}$. The space $H'_q$
is a complement of $B_q$ inside $Z_q$ (so that $H'_q \cong H_q$).
Note that such decompositions depend on the parameter $e$ -- this
dependence is suppressed in the notation here. Given the sequence
of truncation operators $Q^r=(Q^r_d, \cdots ,Q^r_0)$ we use the
decompositions ($\oplus$) to define \emph{modified truncation
operators} $\overline{Q}_q^r$ on the $C_q$'s. In contrast to the
truncation operators, the modified truncation operators are
constructed in such a way as to commute with the boundary maps:
${\partial}_q \circ \overline{Q}^r_q = \overline{Q}^r_{q-1} \circ
{\partial}_q$. Thus $\overline{Q}^r:=(\overline{Q}^r_d, \cdots
,\overline{Q}^r_0)$ is an \emph{endomorphism} of the complex
\begin{gather}
0 \longrightarrow C_d \overset{{\partial}_d }\longrightarrow
\cdots \overset{{\partial}_1}\longrightarrow C_0 \longrightarrow 0
\tag{$\mathcal{C}$}
\end{gather}

For this construction we first use the isomorphism
$B'_{q-1}\overset{\partial}{\underset{\sim}{\longrightarrow}}
B_{q-1}$ to define operators $\widetilde{Q}_{q-1}^r$ on $B'_{q-1}$
in such a way that the following diagram commutes:
\[
\begin{array}{ccc}
    B'_{q-1}         & \overset{\partial}{\underset{\sim}{\longrightarrow}}  & B_{q-1}\qquad \qquad \qquad \\

    \Big\downarrow {\scriptstyle \widetilde{Q}^r_{q-1}}  &                 & \Big\downarrow {\scriptstyle Q_{q-1}^r [B_{q-1};B_{q-1}]} \\

    B'_{q-1}         &
    \overset{\partial}{\underset{\sim}{\longrightarrow}} &
    B_{q-1}\qquad \qquad \qquad
\end{array}
\]
Now if we represent $Q_q^r$ on $ C_q = B_q \oplus H'_q \oplus
B'_{q-1} $ by the matrix:
\[ Q_q^r =  \left( \begin{array}{ccc}
      Q_q^r[B_q;B_q]      &      Q_q^r[B_q;H'_q] & Q_q^r[B_q;B'_{q-1}] \\
      Q_q^r[H'_q;B_q]     &     Q_q^r[H'_q;H'_q] & Q_q^r[H'_q;B'_{q-1}] \\
      Q_q^r[B'_{q-1};B_q] & Q_q^r[B'_{q-1};H'_q] & Q_q^r[B'_{q-1};B'_{q-1}]
          \end{array} \right)            \]
then we define $\overline{Q}_q^r$ by the matrix:
\[ \overline{Q}_q^r := \left( \begin{array}{ccc}
       Q_q^r[B_q;B_q]      &      Q_q^r[B_q;H'_q] & Q_q^r[B_q;B'_{q-1}]\\
        0   &   Q_q^r[H'_q;H'_q] & Q_q^r[H'_q;B'_{q-1}]\\
        0   &   0   & \widetilde{Q}^r_{q-1}
          \end{array} \right)            \]

\begin{claim}
The modified truncation operators commute with the boundary
operators: ${\partial}_q \circ \overline{Q}^r_q =
\overline{Q}^r_{q-1} \circ {\partial}_q$ for $d \geq q \geq 1$.
\end{claim}
\begin{proof}
Write $\omega \in C_q=Z_q \oplus B'_{q-1}$ as $\omega = z + b'$,
where $z\in Z_q$ and $b' \in B'_{q-1}$. Since $\overline{Q}_q^r z
\in Z_q$, applying $\partial$ to it gives $\partial
\overline{Q}_q^r z =0$. Hence we obtain:
\begin{eqnarray*}
\partial \overline{Q}_q^r \omega = \partial \overline{Q}_q^r z +
\partial \overline{Q}_q^r b' & = & \partial \overline{Q}_q^r b'\\ & = &
\partial \widetilde{Q}^r_{q-1} b'\\ & = & Q_{q-1}^r [B_{q-1};B_{q-1}]
\partial b'\\ & = & \overline{Q}_{q-1}^r \partial b' =
\overline{Q}_{q-1}^r \partial \omega
\end{eqnarray*}
\end{proof}
\begin{claim}
The modified truncation operators $\overline{Q}_q^r$ have finite
rank.
\end{claim}
\begin{proof}
Since the original truncation operators $Q_q^r$ have finite rank,
the operators $Q_q^r[*,\star]$ all have finite rank (since they
are all of the form $Q_q^r[*,\star]=P_\star Q_q^r P_*$, where
$P_*$ and $P_\star$ are projection operators). Note that since
$\widetilde{Q}^r_{q-1}$ acts the same (under the above
isomorphism) as $Q_{q-1}^r [B_{q-1};B_{q-1}]$, it also has finite
rank. Now, since all the operator entries in the matrix defining
$\overline{Q}_q^r$ have finite rank, this modified truncation
operator also has finite rank.
\end{proof}
\begin{claim}\label{tendsto}
The modified truncation operators $\overline{Q}_q^r \mbox{ tend to
} {Id}_{C_q}$, as $r$ tends to $\infty$, in the following sense:
$$\forall \omega \in C_q \quad \exists r_{\omega}\in \mathbb{R} \quad
\mbox{s.t.}\quad \overline{Q}^r_q (\omega) = \omega \quad \forall
r \geq r_{\omega}.$$
\end{claim}
\begin{proof}
Given $\omega\in C_q$, write its components with respect to $C_q =
B_q \oplus H'_q \oplus B'_{q-1}$ as $\omega = \omega_1 + \omega_2
+ \omega_3$. Let $r_\omega$ be large enough so that
$support(\omega_i)\subset B(o,r_\omega)$, for $1 \leq i \leq 3$.
We have
\[{\scriptscriptstyle
\begin{array}{ccccccccc}
\omega_1 &=& Q_q^r \omega_1 &=& Q_q^r[B_q;B_q]\omega_1 &+&
Q_q^r[H'_q;B_q]\omega_1 &+& Q_q^r[B'_{q-1};B_q]\omega_1,\\
\omega_2 &=& Q_q^r \omega_2 &=& Q_q^r[B_q;H'_q]\omega_2 &+&
Q_q^r[H'_q;H'_q]\omega_2 &+& Q_q^r[B'_{q-1};H'_q]\omega_2,\\
\omega_3 &=& Q_q^r \omega_3 &=& Q_q^r[B_q;B'_{q-1}]\omega_3 &+&
Q_q^r[H'_q;B'_{q-1}]\omega_3 &+& Q_q^r[B'_{q-1};B'_{q-1}]\omega_3.
\end{array} }\]
Since the decomposition above is a direct sum decomposition, we
see that
\[{\scriptscriptstyle
\begin{array}{ccccccccc}
\omega_1 &=& Q_q^r \omega_1 &=& Q_q^r[B_q;B_q]\omega_1 &+&
0 &+& 0,\\
\omega_2 &=& Q_q^r \omega_2 &=& 0 &+&
Q_q^r[H'_q;H'_q]\omega_2 &+& 0,\\
\omega_3 &=& Q_q^r \omega_3 &=& 0 &+& 0 &+&
Q_q^r[B'_{q-1};B'_{q-1}]\omega_3.
\end{array} }\]
Thus by the definition of $\overline{Q}_q^r$ we have
\[{\scriptscriptstyle
\begin{array}{ccccccccc}
\overline{Q}_q^r \omega_1 &=& Q_q^r[B_q;B_q]\omega_1
&=&\omega_1,\\
\overline{Q}_q^r \omega_2 &=& Q_q^r[H'_q;H'_q]\omega_2
&=&\omega_2,\\
\overline{Q}_q^r \omega_3 &=& \widetilde{Q}_q^r \omega_3. &&
\end{array} }\]
 Note
that since $support(\omega_3)\subset B(o,r_\omega)$ also
$support(\partial \omega_3)\subset B(o,r_\omega)$. Hence
$\partial \widetilde{Q}_q^r \omega_3 = Q_{q-1}^r
\partial \omega_3 = \partial \omega_3$, and since
$B'_{q-1}\overset{\partial}{\underset{\sim}{\longrightarrow}}
B_{q-1}$, we obtain $\overline{Q}_q^r \omega_3=\widetilde{Q}_q^r
\omega_3=\omega_3$. We see that for $1 \leq i \leq 3$,
$\overline{Q}_q^r \omega_i = \omega_i$. That is $\overline{Q}_q^r
\omega = \omega$.
\end{proof}
We summarize the main properties of the modified truncation
operators $\overline{Q}_q^r$:
\begin{itemize}
\item{$\overline{Q}_q^r$ is a finite rank operator.}
\item{The $\overline{Q}_q^r$'s \emph{commute} with the boundary
maps $\partial$.}
\item{$\overline{Q}_q^r \mbox{ tends to } {Id}_{C_q}$, as $r$ tends to
$\infty$.}
\end{itemize}
Hence $\overline{Q}^r$ is an endomorphism of finite rank of the
complex $(\mathcal{C})$.
\begin{rem} Since the complex $(*)$ is exact, the complex $(\mathcal{C})$
 is exact at $C_q$, $d \geq q \geq 1$, and so:
\begin{eqnarray*}
C_0 & = & B_0 \oplus H'_0 \\
C_q & = & B_q \oplus B'_{q-1}, \;\;\; d \geq q \geq 1.
\end{eqnarray*}
\end{rem}

\section{A character formula for functions}
Let $f \in C^{\infty}_c(G)$ be a locally constant function of
compact support on $G$.\\ Recall that the vector spaces $C_q$ are
\emph{smooth} $G$-modules, and let $T_g$ denote the action of $g$
on the $C_q$'s. As is usual for smooth representations, define the
operators $T_f={(T_f)}_q$ on the $C_q$'s:
\[ T_f(\omega) := \underset{G}\int f(g) T_g(\omega)\, dg, \;\;\; \omega\in C_q\]
Note that $T_f=((T_f)_d, \cdots ,(T_f)_0)$ is an endomorphism of
the complex $(\mathcal{C})$ and that the operators ${(T_f)}_q$ are
not necessarily of finite rank: the representation $C_q$ is smooth
but not necessarily admissible.

Since $\overline{Q}^r$ is a finite rank endomorphism of the
complex $\mathcal{C}= (C_d,\cdots,C_0)$, the composition $T_f
\overline{Q}^r=(T_f  \overline{Q}^r_d, \cdots, T_f
\overline{Q}^r_0)$ is also finite rank endomorphism. Recall that
any endomorphism of a complex induces operators on its homology
modules. We use the same notation $T_f  \overline{Q}^r_q$ for the
induced operators on $H_q := H_q(\mathcal{C})$.
\begin{thm}\label{thmfunction} For $G$ connected reductive,
$V$ finitely generated admissible, $e \geq e_0$, and $f \in
C^{\infty}_c(G)$ there exists a radius $r_0(f)$ large enough
\footnote{Choose $\omega_1,\cdots, \omega_n \in H'_0\subset C_0$
such that $\{\omega_1 B_0,\cdots, \omega_n B_0\}\subset H_0 =
C_0/B_0$ is a basis for $W:=T_f (H_0)$. Taking `$r$ large enough'
means that $\omega_i\subset C_c^{or}(X_{(0)}^r;\gamma_e(V))$, for
all $0 \leq i \leq n$.} so that for all $r\geq r_0(f)$:
\[ \sum_{q=0}^{d}  (-1)^q
trace(T_{f}\overline{Q}_q^r ; C_q)=trace(\pi(f) ; V) \]
\end{thm}
\begin{proof}
\begin{eqnarray*}
  \sum_{q=0}^{d}  (-1)^q trace(T_{f}\overline{Q}_q^r ; C_q) & = &
      \sum_{q=0}^{d}  (-1)^q trace(T_{f}\overline{Q}_q^r ; H_q)\\
      & = & trace(T_{f}\overline{Q}_0^r ; H_0)\\
      & = & trace(T_f ; H_0)\\
      & = &  trace(\pi(f) ; V)
\end{eqnarray*}
Applying the Hopf trace formula for finite rank operators
~\cite[Proposition 2.1]{AB1} to the finite rank endomorphism
$T_{f}\overline{Q}^r$ we obtain the first equality. Since the
complex ($\mathcal{C}$) is exact at all places $C_q,\; d \geq q
\geq 1$, the only (possible) non-zero homology is $H_0$. This
explains the second equality. Using the exact $G$-equivariant
complex $(*)$ we have that $H_0 = C_0/ B_0 \cong V$ and hence that
$trace(T_f ; H_0) = trace(\pi(f) ; V)$. Note that since $V$ is
admissible the operator $\pi(f)$ is of finite rank and so taking
its trace makes sense. This gives the fourth equality. The
following explains the
third equality.\\

Consider the diagram:
\[
\begin{array}{ccc}
    H_0  & \overset{\epsilon}{\underset{\sim}{\longrightarrow}} & V\\

    \Big\downarrow {T_f}  &   & \Big\downarrow {\pi(f)}\\

    H_0  & \overset{\epsilon}{\underset{\sim}{\longrightarrow}} & V\\
\end{array}
\]
Since the exact complex $(*)$ respects the action of $g$, for all
$g\in G$, it also respects the action of $f$; therefore the
diagram above commutes. Since $V$ is admissible, $\pi(f)$ is a
finite rank operator on $V$. Hence $T_f$ is also a finite rank
operator on $H_0$. To show
\[ trace\, (T_f \overline{Q}^r_0;H_0)  = trace\, (T_f ;H_0), \]
it is enough to show that for large enough $r$, the operators
$\overline{Q}^r_0$ act trivially on the finite dimensional image
$T_f (H_0)$ of $T_f$.\\
Choose $\omega_1,\cdots, \omega_n \in H'_0\subset C_0$ such that
$\{\omega_1 B_0,\cdots, \omega_n B_0\}\subset H_0 = C_0/B_0$ is a
basis for $W:=T_f (H_0)$, and let $W'$ be a complement of $W$
inside $H_0$: $H_0=W' \oplus W.$ With respect to this direct sum
decomposition the operator $T_f$ is represented by a matrix of the
form:
\[ T_f = \left( \begin{array}{cc}
      0 & 0 \\
      * & T_{W,W} \\
          \end{array} \right).            \]
Let $r$ be large enough so that $\omega_i\subset
C_c^{or}(X_{(0)}^r;\gamma_e(V))$, for all $0 \leq i \leq n$. Then
(Claim ~\ref{tendsto}) $\overline{Q}_0^r(\omega_i) = \omega_i $.
Thus $\overline{Q}_0^r(\omega_i B_0)=\overline{Q}_0^r(\omega_i)B_0
= \omega_i B_0$, and so
$${\overline{Q}^r_0|}_W \equiv {Id|}_W.$$ Note that $r$ depends on
$W$ which depends on $f$.\\
So with respect to the above decomposition of $H_0$ the operator
${\overline{Q}^r_0}$ is represented by a matrix of the form:
\[ \overline{Q}^r_0 = \left( \begin{array}{cc}
      * & 0 \\
      * & 1 \\
          \end{array} \right).            \]
Multiply $T_f$ and $\overline{Q}^r_0$ to obtain:
\[ T_f \overline{Q}^r_0 = \left( \begin{array}{cc}
      0 & 0 \\
      * & T_{W,W} \\
          \end{array} \right)
   \left( \begin{array}{cc}
      * & 0 \\
      * & 1 \\
          \end{array} \right) =    \left( \begin{array}{cc}
      0 & 0 \\
      * & T_{W,W} \\
          \end{array} \right).          \]
Equality of the traces follows:
\[ trace\, T_f \overline{Q}^r_0  = trace\, T_f \]
This completes the proof of the theorem.
\end{proof}

\subsection*{How to proceed} Keeping in mind the formula of Theorem
~\ref{thmfunction}, we have the following two objectives in mind:
\begin{enumerate}
\item{} To give a more geometric interpretation of this formula in terms
of the original truncation operators.
\item{} To replace the function $f$ in the formula with a compact
element $\gamma$.
\end{enumerate}
In the rest of this paper we show how to complete step $(1)$ for
groups of semisimple rank $1$; we make some progress in the
direction of showing $(1)$ for a general reductive group. Assuming
step $(1)$ we show how to complete step $(2)$ for all compact
elements. For groups of semisimple rank $1$, we recover the
Schneider-Stuhler result.

\section{Nice decompositions}\label{nicedecomp}
Let $\gamma\in G^{cpt}$ and fix $o\in X^\gamma$. Recall the
decomposition: $C_q=C_q^{in(r)} \oplus C_q^{out(r)}$ which is
$T_\gamma$-equivariant (see paragraph after Example ~\ref{SL2ex}).
Let $B_q^{in(r)}:=B_q \cap C_q^{in(r)}$. It is a $T_\gamma$-stable
subspace of $B_q$ as $T_\gamma$ acts on $B_q$ and on
$C_q^{in(r)}$. Let $B_{q}^{out(r)}$ be any $T_\gamma$-stable
complement of $B_{q}^{in(r)}$ in $B_{q}$
(we will show existence of such complements later in this section).\\

\begin{defn}\label{def-nice-decomp} We say that the decomposition ($\oplus$) is \emph{`nice'}
relative to $\gamma\in G^{cpt}$ and $r\in \mathbb{R}$, (here
$\gamma$ is such that $o \in X^\gamma$), if it has the following
form:
\begin{itemize}

\item{$C_0  =B_0 \oplus H'_0$}

\item{$C_q  =B_q \oplus
B'_{q-1},  \quad  d \geq q \geq 1$.\\
Where $B'_{q-1} = (B_{q-1}^{in(r)})' \oplus (B_{q-1}^{out(r)})'$,
for some subspaces $(B_{q-1}^{in(r)})'$ and $(B_{q-1}^{out(r)})'$
which satisfy:\\ $(B_{q-1}^{in(r)})'
\overset{\partial}{\underset{\sim}{\longrightarrow}}(B_{q-1}^{in(r)}),
\;\;(B_{q-1}^{out(r)})'
\overset{\partial}{\underset{\sim}{\longrightarrow}}
(B_{q-1}^{out(r)})$, and\\ $(B_{q-1}^{in(r)})'\subset
C_q^{in(r)}$.}

\item{All $\oplus$ decompositions here are $T_\gamma$-equivariant.}
\end{itemize}

In this definition the dependence on the parameter $e$ is
suppressed in the notation. If for each $r\in\mathbb{R}$ nice
decompositions exist for all $e \geq e_0$, we will say that the
nice decomposition is \textit{uniform in $e$}. If for each $r$
nice decompositions exist only for all $e \geq e_r$, where $e_r$
is an integer depending in $r$, we will say that the nice
decomposition \textit{depends on $e_r$}.

\end{defn}
In this section we show that for groups of semisimple rank $1$,
such nice decompositions exist uniformly in $e$. For a general
connected reductive group we show that nice decompositions exist,
but that they depend on $e_r$.\\

When the decomposition $(\oplus)$ is nice, we will show, under
appropriate conditions, that:
$$trace(T_\gamma Q_q^r; C_q) = trace(T_\gamma \overline{Q}_q^r; C_q)
\quad  d \geq q \geq 0$$

\subsection{Truncated complexes}\label{trunc-compl} For a general group $G$ the direct
sum decomposition is controlled by the following \emph{truncated
complex}:
\begin{gather}
 0 \longrightarrow C_c^{or}(X_{(d)}^r ;\gamma_e(V)) \overset{\partial}\longrightarrow
 \cdots
\overset{\partial}\longrightarrow C_c^{or}(X_{(0)}^r ;\gamma_e(V))
\overset{\epsilon}\longrightarrow V   \tag{$\overline{*}$}
\end{gather}

The exactness of this complex will guarantee the existence of a
nice direct sum decomposition $(\oplus)$. The exactness will
follow from the properties of the truncated building $X^r$.
Ideally we would like to show exactness of $(\overline{*})$ for
all $e \geq e_0$, but the techniques used here will show (for a
general connected reductive group) exactness
only for $e \geq e_r$, for some $e_r$ which depends on $r$.\\

We recall the following averaging process:\\
Let $K \subset G$ be an open compact subgroup and denote by
$1_{K}$ its characteristic function. The operator $\pi (K):=
\pi(1_{K})$ on $V$ is well defined as $1_{K} \in C^{\infty}_c(G)$:
\[\pi(K)v = \underset{k\in K}\int \pi(k)v\, dk, \;\;\; v\in V\]
where $dk$ is a normalized Haar measure on $K$.\\

Note the following facts about $\pi(K)v$:
\begin{itemize}

\item{For any $v\in V$, $\pi(K)v \in V^K$.\\
This follows from the definition of $\pi(K)$.}

\item{If for some open compact subgroup $U$, $v\in V^U$ and $K$ normalizes
$U$, then $\pi(K)v \in V^U$:
   \begin{eqnarray*}
      \pi(u) \pi(K)v  & = & \pi(u) \underset{k\in K}\int \pi(k)v\,
      dk
        =  \underset{k\in K}\int \pi(u) \pi(k)v\, dk\\
      & = & \underset{k\in K}\int \pi(k) \pi(k^{-1}) \pi(u) \pi(k)v\, dk
        =  \underset{k\in K}\int \pi(k) \pi(k^{-1}uk)v\, dk\\
      & = & \underset{k\in K}\int \pi(k) \pi(u')v\, dk
        =   \underset{k\in K}\int \pi(k) v\, dk
        =  \pi(K)v
  \end{eqnarray*}
     }
\end{itemize}

Recall that for any facet $F$, the groups $U_F^{(e)}$ are
normal subgroups of  $P^{\dagger}_F$. For vertices $x$ this means
that $U_x^{(e)} \lhd P_x = P^{\dagger}_x$.

Let $X^r$ be the truncated building with center $o \in X$.
\begin{claim} It is possible to choose an integer $e_r$ large enough so that
for any $e \geq e_r$ the groups $U_x^{(e)}$, for $x \in X^r$, all
normalize each other.
\end{claim}
\begin{proof}
Let $U(r):= \underset{y\in X^r}\bigcap U^{(0)}_y$, an open compact
subgroup (since this intersection is a finite intersection of open
compact subgroups). Since for any vertex $x$, the subgroups
$U_x^{(e)}$ form a filtration of the identity element in $G$, we
can choose $e_r$ large enough so that $U_x^{(e)} \subset U(r)
\subset P_y$, for all $x,y \in X^r$ and $e \geq e_r$. Since
$U^{(e)}_x \subset P_y$, it normalizes $U^{(e)}_y$. As this holds
for all $x,y \in X^r$ and $e \geq e_r$, the claim is proved.
\end{proof}

\begin{thm} For $e \geq e_r$, the sequence $(\overline{*})$ is exact.
\end{thm}
\begin{proof}
Take a non-zero cycle $\omega \in C_c^{or}(X_{(q)}^r
;\gamma_e(V))$; we need to show that it is a boundary. The proof
is by induction
on the number of facets in the support of $\omega$.\\

We can use any open compact subgroup, $K \subset U(r)$, to average
$\omega$:
\[  (\pi(K)\omega )((F',c'))= \underset{k\in K}\int \pi(k)(\omega ((k^{-1}F',k^{-1}c')))\, dk \]
Since $K \subset P_{F'}$ for all $F' \subset X^r$, $k\in K$ acts
trivially on $(F',c') \in X^r_q$. Hence
\[  (\pi(K)\omega )((F',c'))= \underset{k\in K}\int \pi(k)(\omega ((F',c')))\, dk .\]

Note that $ (\pi(K)\omega )((F',c')) \in V^K \cap
V^{U^{(e)}_{F'}}$ and that $\pi(K)\omega \neq 0$. Also, since the
action of $K\subset G$
commutes with the boundary operators, $\pi(K)\omega$ is also a cycle.  \\

Let $F\in X^r_{(q)}$ be an oriented $q$-facet such that $F \subset
support(\omega)$. Use $K=U^{(e)}_F$. Then $\omega =
\pi(U^{(e)}_F)\omega + (\omega - \pi(U^{(e)}_F)\omega)$. Being a
difference of two cycles, $(\omega - \pi(U^{(e)}_F)\omega)$ is
itself a cycle. Note that $\omega(F)\in V^{U_F^{(e)}}$ implies
that $(\pi(U_F^{(e)})\omega)(F)=\omega(F)$, hence $(\omega -
\pi(U_F^{(e)})\omega)(F)=0$. This means that the support of
$\omega - \pi(U_F^{(e)})\omega$ is strictly smaller than the
support of $\omega$. Thus by induction we conclude that $\omega -
\pi(U_F^{(e)})\omega$ is a boundary. The non-zero cycle
$\pi(U^{(e)}_F)\omega$ has the same support as $\omega$ and is
invariant (by the averaging process) under $U^{(e)}_F$. Since $e
\geq e_r$, the groups $U_F^{(e)}$, $F\in X^r_{(q)}$, all normalize
each other, and so $\pi(U_{F'}^{(e)})\pi(U_F^{(e)})\omega$ is
still $U_F^{(e)}$-invariant. Thus after applying this process a
finite number of times, i.e. as $K$ ranges over all $F \subset
support(\omega)$, we can assume with out loss of generality that
the cycle $\omega$ is invariant under all $U_F^{(e)}$, $F \subset
support(\omega)$. \\

Denote by $X_0(\omega)$ the set of vertices in
$\overline{support(\omega)}$ and let $$V^{(e)}_{\omega}:=
\underset{x \in X_0(\omega)}\bigcap V^{U^{(e)}_x},$$ then $\omega
((F',c')) \in V^{(e)}_{\omega}$, all $(F',c') \subset X^r$.\\


We show that this situation can be reduced to a constant
coefficients case.\\

Let $S \subset X$ be any subset of the building. We recall the
$simplicial$ operation and define an algorithmic version of the
$convex$ operation:

\begin{enumerate}
\item{}{ $geod(S):= \{z \in X \;|\; z \in geod(x,y) \mbox{ for some } x,y \in
S\}$.}\\
\item{}{ $simp(S):= \emph{smallest subcomplex of the building containing S}$.}
\end{enumerate}
\begin{lem} If $V^{(e)}_{\omega} \subset V^{U_z^{(e)}}, \mbox{ for all
} z \in S$ then
\begin{enumerate}
\item{} $V^{(e)}_{\omega} \subset V^{U_z^{(e)}}, \mbox{ for
all } z \in geod(S)$.
\item{} $V^{(e)}_{\omega} \subset V^{U_z^{(e)}}, \mbox{ for all }
z \in simp(S)$.
\end{enumerate}
\end{lem}
\begin{proof}
For $z \in geod(x,y)$ have $U_z^{(e)} \subset U_x^{(e)}U_y^{(e)} $
~\cite[Lemma 1.28]{Vig}. Hence $V^{U_z^{(e)}} \supset
V^{U_x^{(e)}} \cap V^{U_y^{(e)}} \supset V^{(e)}_{\omega}$, which
proves $(1)$. If $z \in simp(S)$, then $z \in \overline{F}$ for
some facet $F$ s.t $F \cap S \neq \emptyset$. So (by property
$(U4)$) $U_z^{(e)} \subset U^{(e)}_F = U^{(e)}_x$ for any $x \in F
\cap S$ and hence $V^{(e)}_{\omega} \subset V^{U_x^{(e)}} \subset
V^{U_z^{(e)}}$, which proves $(2)$.
\end{proof}

We continue with the proof of the theorem, where $X^r$ is the
truncated building with parameter $r$ and we choose $e_r$ as in
the claim above. Let $\omega \in C_c^{or}(X_{(q)}^r
;\gamma_{e_r}(V))$ be a $q-$cycle. We want to show it is a
boundary. Let $S_{\omega} := supp(\omega)$ and consider the
sequence of increasing subsets:
\begin{gather}
 S_{\omega} \subseteq geod(S_{\omega}) \subseteq simp(geod(S_{\omega})) \subseteq
\cdots \tag{$\ddagger$}
\end{gather}

Essentially the same argument as that showing that the sequence
$(\dagger)$ stabilizes, shows that this sequence $(\ddagger)$ also
stabilizes after a finite number of terms. Denote by
$\overline{S}_\omega$ the stable terms in the sequence
$(\ddagger)$.

It follows from the construction of $\overline{S}_\omega$ that
$\overline{S}_\omega \subset X^r$, and that $\overline{S}_\omega$
is simplicial and convex.

Also since $V^{(e)}_{\omega} \subset V^{U_x^{(e)}}, \mbox{ for all
} x \in S_\omega$, we can apply the lemma several times to
conclude that $V^{(e)}_{\omega} \subset V^{U_x^{(e)}}, \mbox{ for
all } x \in \overline{S}_\omega$. Now consider the following
commutative diagram of chain complexes:

\[
\begin{array}{ccccccc}
     \overset{\partial}\rightarrow
     C_c^{or}(X^r_{(q+1)} ;\gamma_{e}(V)) & \overset{\partial}\rightarrow  &
     C_c^{or}(X^r_{(q)} ;\gamma_{e}(V))   & \overset{\partial}\rightarrow & C_c^{or}(X^r_{(q-1)}
     ;\gamma_{e}(V)) & \overset{\partial}\rightarrow
     \\ \\

     \bigcup&   & \bigcup &    & \bigcup &\\ \\

     \overset{\partial}\rightarrow
     C_c^{or}((\overline{S}_\omega)^r_{(q+1)} ;V^{(e)}_\omega) & \overset{\partial}\rightarrow  &
     C_c^{or}((\overline{S}_\omega)^r_{(q)} ;V^{(e)}_\omega)   & \overset{\partial}\rightarrow &
     C_c^{or}((\overline{S}_\omega)^r_{(q-1)} ;V^{(e)}_\omega) &  \overset{\partial}\rightarrow
     \\ \\

\end{array}
\]

The bottom line is the homology chain complex of the simplicial
complex $\overline{S}_\omega$ with constant coefficients
$V^{(e)}_\omega$. Since $\overline{S}_\omega$ is convex (hence
contractible) its constant coefficient chain complex is exact.
Hence $\omega$, considered as a cycle in
$C_c^{or}((\overline{S}_\omega)^r_{(q)} ;V^{(e)}_\omega))$ must be
a boundary of some chain $\delta \in
C_c^{or}((\overline{S}_\omega)^r_{(q+1)} ;V^{(e)}_\omega)) \subset
C_c^{or}(X^r_{(q+1)} ;\gamma_{e_r}(V)) $.\\ We have found $\delta
\in C_c^{or}(X^r_{(q+1)} ;\gamma_{e_r}(V))$ such that $\partial
\delta = \omega$, i.e., the cycle $\omega$ is a boundary. This
concludes the proof of the theorem.
\end{proof}

\subsection{Existence of nice decompositions}

Let $\gamma$ be a compact element such that $o \in X^\gamma$. The
compact element $\gamma$ is contained in some compact subgroup $K$
of $G$. Let $K_\gamma:=\overline{<\gamma>}$ be the closure in $K$
of the subgroup generated by $\gamma$. $K_\gamma$ is a compact
subgroup containing $\gamma$. It will be used for `averaging'
purposes.
\begin{rem} Schneider-Stuhler make the assumption that
$Z^{\circ}(G)$ is anisotropic which implies that the stabilizer
$P_x$, of any vertex $x$, is a compact group. In particular $P_o$
being a compact group containing $\gamma$ could also be used for
`averaging' purposes. We will not need it here, so we will not
need to make the above assumption on $G$.
\end{rem}
Recall that for $B_{q}^{out(r)}$ we needed to take any
$T_\gamma$-stable complement of $B_{q}^{in(r)}$ in $B_{q}$.\\

To see such a complement exists, take any projection $\pi :B_q
\rightarrow B_{q}^{in(r)}$ and average it over the compact group
$K_\gamma$:

$$\int\limits_{K_\gamma} T_k^{-1} \circ \pi \circ T_k$$

This makes sense since $K_\gamma$ acts on $B_q$ and on
$B_{q}^{in(r)}$. The result is a $T_\gamma$-equivariant projection
onto $B_{q}^{in(r)}$.
Let $B_{q}^{out(r)}$ be the kernel of this projection.\\

Get a $T_\gamma$-equivariant direct sum decomposition of $B_q$ for
each $d \geq q \geq 0$:  $$B_q=B_q^{in(r)} \oplus B_q^{out(r)}$$

Note that $B_{q}^{in(r)} \subset C_{q}^{in(r)}$ but
$B_{q}^{out(r)} \not \subset C_{q}^{out(r)}$.

Recall that since the chain complex $\mathcal{C}$ is exact at
$C_q$, $d \geq q \geq 1$, the direct sum decomposition ($\oplus$)
has the form :
\begin{eqnarray*}
C_0 & = & B_0 \oplus H'_0 \\
C_q & = & B_q \oplus B'_{q-1}, \;\;\; d \geq q \geq 1.
\end{eqnarray*}

Here $H_0'$ is any $T_\gamma$-equivariant complement of $B_0$ in
$C_0$. As above, such a complement can be realized as the kernel
of an appropriate $T_\gamma$-equivariant projection from
$C_0$ to $B_0$ which can be produced by averaging over $K_\gamma$.\\
For $B'_{q-1}$ we will chose a particularly `nice' section of
$\partial$ in the following sense.

 For $q \geq 1$ we will construct a particular section
$\alpha$, i.e. $\partial \circ \alpha = id$, of the surjective map
$\partial:C_{q} \longrightarrow B_{q-1}$ and then use it to define
$B_{q-1}'$ as the image $\alpha(B_{q-1})$.

\begin{ex} Let $G=SL_2$ and $V=\mathbb{C}$ be the trivial
representation. The building $X$ is a tree. The chain complex
$(\mathcal{C}) $ for $X$ is: $0 \longrightarrow C_1
\longrightarrow C_0 \longrightarrow 0$.\\ Let $o$ be a vertex
fixed by $\gamma$. A convenient basis for $B_0$ is $\{\delta _x -
\delta _o \}$ where $x$ runs over all vertices $\neq o$.\\
Define the section $\alpha$ of $\partial$ to be $\alpha (\delta _x
- \delta _o)=$ sum of all the edges connecting $o$ to $x$.

Note that $\alpha$ commutes with the action of $\gamma$ so that it
is $T_\gamma$-equivariant, and also that $\alpha (\delta _x -
\delta _o)\in C_q^{in(r)}$ if $x\in X^r$.
\end{ex}

The point of the following lemma is to show that in general we can
always find a section with the well behaved properties $\alpha$
has in this example.
\begin{lem}\label{key lemma} (Key Lemma)

\begin{description}
\item[semisimple rank $1$ case] Suppose $G$ has semisimple rank $1$. Suppose $\gamma \in G$ is such that
$o \in X^\gamma$. For each integer $r \geq 0$ and for all $e \geq
e_0$ the map $\partial :C_{1} \longrightarrow B_{0}$ has a
$T_\gamma$-equivariant section $\alpha$, such that
$\alpha(B_{0}^{in(r)})\subset C_{1}^{in(r)}$.\\

\item[general case] Let $G$ be a connected reductive
group. Suppose $\gamma \in G$ is such that $o \in X^\gamma$. For
each integer $r \geq 0$ there exists an integer $e_r$ so that for
all $e \geq e_r$ each map $\partial :C_{q} \longrightarrow
B_{q-1}$, $d \geq q \geq 1$, has a $T_\gamma$-equivariant section
$\alpha$, such that $\alpha(B_{q-1}^{in(r)})\subset
C_{q}^{in(r)}$.
\end{description}
\end{lem}
\begin{rem} The main difference between the general version of the
key lemma and the semisimple rank $1$ version is the dependence of
the parameter $e$ on the radius of truncation $r$. This dependence
is what makes the general version (significantly?) weaker and
prevents us from achieving steps $(1)$ and $(2)$ in general.
\end{rem}
\begin{proof} $\;C_{q}^{in(r)} \subset C_{q}$ hence $\partial
C_{q}^{in(r)} \subset \partial C_{q} = B_{q-1}$. Since
$C_{q}^{in(r)}$ consists of maps supported on $X^{r}$ their
boundary is also supported on
$X^{r}$, so $\partial C_{q}^{in(r)} \subset C_{q-1}^{in(r)}$.\\
Put together we see that $\partial C_{q}^{in(r)} \subset B_{q-1}
\cap C_{q-1}^{in(r)} =: B_{q-1}^{in(r)}$. So the map $\partial
:C_{q} \longrightarrow B_{q-1}$ restricts to a map $\partial '
:C_{q}^{in(r)} \longrightarrow B_{q-1}^{in(r)}$, which is
$T_\gamma$-equivariant. Recall that $C_{q}^{in(r)} = C_q(X^r)$.
\begin{claim} (semisimple rank $1$ case) Suppose $G$ has semisimple
rank $1$. For each $r \geq 0$ and for each $e \geq e_0$, the map
$\partial'$ is surjective.
\end{claim}

\begin{proof}
Recall the list of facts (section ~\ref{Ugroups}) regarding the
groups
$U_F^{(e)}$.\\
Since in this case the chain complex $(\mathcal{C})$ is: $0
\longrightarrow C_1 \overset{\partial}\longrightarrow C_0
\longrightarrow 0$, we only have to show that $\partial '
:C_{1}^{in(r)} \longrightarrow B_{0}^{in(r)}$ is surjective. That
is, given a non-vanishing $0$-cycle $\omega \in
B_{0}^{in(r)}=C_c^{or}(X_{(0)}^r ;\gamma_e(V)) \cap B_0(X)$, we
need to show that it is a boundary of a chain in
$C_c^{or}(X_{(1)}^r ;\gamma_e(V))$. The proof is by induction on
the radius of the support of $\omega$. Let $B(o,r_{\omega})$ be
the smallest ball containing the support of $\omega$, we will show
that by adding a boundary to $\omega$, we obtain a $0$-cycle
$\omega^{\star}$ with support contained in a smaller ball
$B(o,r_{\omega}-1)$. Repeating this process will show that
$\omega$ is a boundary.
\begin{defn} Given a finite collection of points in the building, $S
\subset X$, we say that a point $x \in S$ is an \emph{extreme
point} (relative to a point $o$), if it is farthest away from $o$.
That is if $d(o,x) \geq d(o,y)$ for all $y \in S$.
\end{defn}
The support of $\omega$ is finite so it makes sense to talk about
an extreme vertex in this support. Let $x$ be such a vertex. Let
$y$ be another vertex in the support of $\omega$ and let $v=
\omega(x) \in V^{U_x^{(e)}}$ and $w= \omega(y) \in V^{U_y^{(e)}}$.
Consider the average of $w$, $\pi(U_x^{(e)})w$, over the compact
group $U_x^{(e)}$.\\

Let $F$ be the facet of the geodesic $geod(y,x)$ whose closure
$\overline{F}$ contains $x$. We next show that $\pi(U_x^{(e)})w
\in V^{U_F^{(e)}}$. For any $u_z \in U_z^{(e)}$, where $z\in F$,
we have:
\begin{eqnarray*}
      \pi(u_z) \pi(U_x^{(e)})w  & = &
                  \pi(u_z) \pi(U_x^{(e)})\pi(u_z^{-1})\pi(u_z)w\\
      & = & \pi(u_z U_x^{(e)} u_z^{-1})\pi(u_z)w\\
      & = & \pi(U_x^{(e)})\pi(u_z)w\\
      & = & \pi(U_x^{(e)})\pi(u_x)\pi(u_y)w\\
      & = & \pi(U_x^{(e)})w.
  \end{eqnarray*}
The first and second equalities are clear. By fact $(U5)$ (section
~\ref{Ugroups}) we see that $u_z$ normalizes $U_x^{(e)}$. This
explains the third equation. Using $(U7)$ we see that $u_z$ is of
the form $u_x u_y$ for some $u_x \in U_x^{(e)}$ and some $u_y \in
U_y^{(e)}$, which gives the forth equality. By absorbing $u_x$
into $U_x^{(e)}$ and since $u_y$ acts trivially an $w$ we obtain
the fifth equality. Thus
$\pi(U_x^{(e)})w \in V^{U_F^{(e)}}$.\\

\begin{figure}[h]
\psfragscanon \centering \psfrag{x}{$x$} \psfrag{x1}{$x_1$}
\psfrag{x2}{$x_2$} \psfrag{y1}{$y_1$} \psfrag{y2}{$y_2$}
\psfrag{F1}{$F_1$} \psfrag{F2}{$F_2$} \psfrag{o}{$o$}
\epsfig{file=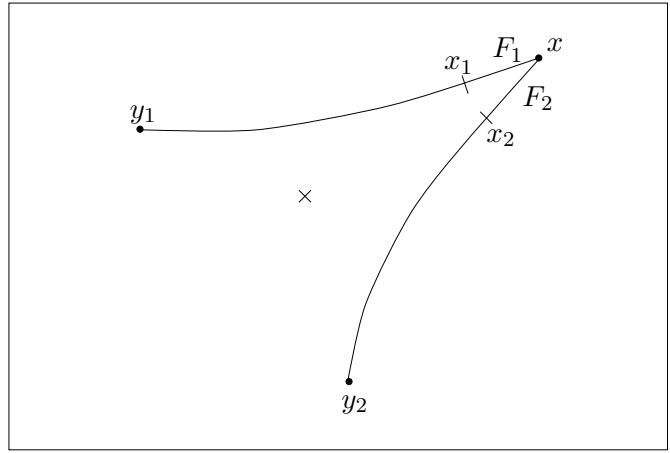,height=6cm}\caption{\label{geodesics}Geodesics
to $x$}
\end{figure}

Now label the vertices in the support of $\omega$ which are
distinct from $x$ by $y_1, \cdots ,y_n$, (see Figure
~\ref{geodesics}), and let $w_i=\omega(y_i)$. Let $F_i$ denote the
facet of the geodesic $geod(y_i,x)$ whose closure $\overline{F_i}$
contains $x$, oriented as to point towards $x$. Since $\omega$ is
a $0$-cycle, we have:
$$\epsilon(\omega)=v + w_1 + \cdots + w_n = 0,$$
where $\epsilon$ is the augmentation map of the complex $(*)$.\\
Applying $\pi(U_x^{(e)})$ to each term above we obtain:
$$v + \pi(U_x^{(e)})(w_1) + \cdots + \pi(U_x^{(e)})(w_n) = 0.$$
Let $\overline{w}_i := \pi(U_x^{(e)})(w_i)$. Since
$\pi(U_x^{(e)})(w_i) \in V^{U_{F_i}^{(e)}}$, the $1$-chain:
$\omega_x := \delta_{F_1}^{\overline{w}_1} + \cdots +
\delta_{F_n}^{\overline{w}_n} $ is in $C_{1}^{in(r)}$ (where
$\delta_F^v$ is a `delta chain' supported on $F$ with value $v$).
Now, $\omega = (\omega +
\partial \omega_x) -
\partial \omega_x$, where $\omega + \partial \omega_x \in
B_{0}^{in(r)}$. Note that $(\omega + \partial \omega_x)(x) = 0$,
so that $\omega +
\partial \omega_x$ is not supported on $x$. Note also that as a
consequence of this process, we introduced into the  support of
$\omega +
\partial \omega_x$ the (possibly) new vertices $x_i$, where
$\partial \delta_{F_i}^{\overline{w}_1} =
\delta_{x}^{\overline{w}_1} - \delta_{x_i}^{\overline{w}_1}$. We
next show that the vertices $x_i$ are strictly closer to $o$ than
$x$ and so $\omega +
\partial \omega_x$ is an improvement on $\omega$ in the sense of
trying to show it is a boundary. Now repeating the process above
for all extreme vertices of distance $d(o,x)$ we obtain the
$0$-cycle $\omega^{\star}$ which is supported on
$B(o,r_{\omega}-1)$. Thus by induction we have shown that $\omega$
is a boundary of a $1$-chain in $C_{1}^{in(r)}$ and the claim is
proved.
\begin{rem} In the proof above, since the building is $1$-dimensional,
all the facets $F_i$ are actually the same facet. Since parts of
this proof generalize to higher dimensional buildings we leave it
intact.
\end{rem}
In order to show that this process terminates we need to show that
the vertices $x_i$ are strictly closer to $o$ than $x$.\\ We use
the \emph{negative curvature inequality}, ~\cite[p.153]{Brown},
which says that for all $t \in [0,1]$:
$$d^2(o,p_t) \leq (1-t)d^2(o,y)+td^2(o,x)-t(t-1)d^2(y,x)$$
where $p_t:=(1-t)y+tx \in geod(y,x)$, $t\in [0,1]$. For $t\in
(0,1)$ we have:
$$d^2(o,p_t) \lneq (1-t)d^2(o,y)+td^2(o,x)$$
as $t(t-1)d^2(y,x) \gneq 0$. Now $d(o,y) \leq d(o,x)=r_{\omega}$,
so:
$$d^2(o,p_t) \lneq (1-t)r_{\omega}^2+tr_{\omega}^2=r_{\omega}^2$$
hence:
$$d(o,p_t) \lneq r_{\omega}.$$
Since $x_i=p_t$ for some $t\in(0,1)$ we see that the vertices
$x_i$ are indeed closer to $o$ than is $x$:
$$d(o,x_i) \lneq d(o,x).$$
\end{proof}
\begin{claim}(general case) Let $G$ be a connected reductive
group. For each $r \geq 0$ there exists an integer $e_r$,
depending on $r$, such that for all $e \geq e_r$ the map $\partial
':C_{q}^{in(r)} \longrightarrow B_{q-1}^{in(r)}$ is surjective.
\end{claim}
\begin{proof}
We have the following commutative diagram:
\[
\begin{array}{ccccc}
     C_q(X) & \overset{\partial}\longrightarrow  &
     B_{q-1}(X) & & \\

     \bigcup &   & \bigcup & & \\

     C_q(X^r) & \overset{\partial'}\longrightarrow  &
     B_{q-1}^{in(r)} & := & B_{q-1}(X) \cap C_{q-1}(X^r)\\
\end{array}
\]

The map $\partial: C_q(X) \longrightarrow  B_{q-1}(X)$ is
surjective by definition. Using $B_{q-1}(X^r)  \subset
B_{q-1}^{in} \subset B_{q-1}(X)$ and the exactness of
$(\overline{*})$ we obtain:
$$B_{q-1}(X^r)= Z_{q-1}(X^r)=B_{q-1}(X) \cap Z_{q-1}(X^r)=
B_{q-1}(X) \cap C_{q-1}(X^r) =: B_{q-1}^{in}$$

Thus the map $\partial' :C_q^{in(r)} \longrightarrow
B_{q-1}^{in(r)}$ is surjective and the claim is proved.
\end{proof}
Now, continuing with the proof of the key lemma, let $\alpha ':
B_{q-1}^{in(r)} \longrightarrow C_q^{in(r)}$ be a
$T_\gamma$-equivariant section of $\partial '$ and let $\alpha
'':B_{q-1}^{out(r)}\longrightarrow C_{q}(X)$ be any section of
$\partial:C_{q}(X) \longrightarrow B_{q-1}(X)$. Again use
$K_\gamma$ to average $\alpha ''$ and make it
$T_\gamma$-equivariant. Let $\alpha:B_{q-1}^{in(r)}(X) \oplus
B_{q-1}^{out(r)}(X)\longrightarrow C_{q}(X)$ be $ \alpha
'\oplus\alpha ''$. $\alpha$ satisfies the properties in the
statement of the lemma.
\end{proof}

\begin{cor}\label{directsumexist} Recall Definition ~\ref{def-nice-decomp} of nice
$\oplus$ decompositions.
\begin{description}
\item[semisimple rank $1$ case] For $G$ of semisimple rank $1$,
nice $\oplus$ decompositions exist and are uniform in $e$.
\item[general case] For $G$ a connected reductive group,
nice $\oplus$ decompositions exist but they depend on $e_r$.
\end{description}
\end{cor}
\begin{proof}
Let $H_0$ be any $T_\gamma$-equivariant complement of $B_0$ inside
$C_0$. Recall the $T_\gamma$-equivariant decomposition $B_{q-1}=
B_{q-1}^{in(r)} \oplus B_{q-1}^{out(r)}$. Use the
$T_\gamma$-equivariant section $\alpha$ of Lemma ~\ref{key lemma}
to define:
$$(B_{q-1}^{in(r)})':=\alpha(B_{q-1}^{in(r)})$$
$$(B_{q-1}^{out(r)})':=\alpha(B_{q-1}^{out(r)}).$$
The properties of $\alpha$ guarantee that this construction gives
a nice $\oplus$ decomposition.
\end{proof}

\subsection{Traces with respect to nice decompositions} Given direct sum
decompositions of the vector spaces $C_q$ ($d \geq q \geq 0$):
\begin{gather}
C_q = B_q \oplus H'_q \oplus B'_{q-1} \tag{$\oplus$}
\end{gather}
we now show that if such decompositions are nice (with respect to
$\gamma$ and $r$) then the finite rank operators $T_\gamma Q_q^r$
and $T_\gamma \overline{Q}_q^r$ all have the same trace.
\begin{prop}\label{ridofbar} Let $\gamma \in G$ be an element such that
$o \in X^\gamma$. If the decomposition ($\oplus$) is nice relative
to $\gamma$ and $r$, then
$$trace(T_\gamma Q_q^r; C_q) = trace(T_\gamma \overline{Q}_q^r; C_q)
\quad  d \geq q \geq 0$$
\end{prop}
\begin{proof}

For $q=0$ have $C_0 = B_0 \oplus H'_0$. Since $B_0$ and $H'_0$ are
$T_\gamma$-stable, the linear operator $T_\gamma$ on $C_0 = B_0
\oplus H'_0$ is block diagonal:

\[ T_\gamma = \left( \begin{array}{cc}
      * & 0 \\
      0 & * \\
          \end{array} \right)            \]

As $\overline{Q}_0^r$ differs from $Q_0^r$ only on the bottom left
block (with respect to $C_0 = B_0 \oplus H'_0)$:

\[ \begin{array}{cc}

 Q_0^r = \left( \begin{array}{cc}
      Q_0^r[B_0;B_0] & Q_0^r[B_0;H'_0] \\
      Q_0^r[H'_0;B_0] & Q_0^r[H'_0;H'_0] \\
          \end{array} \right)
 &
 \overline{Q}_0^r = \left( \begin{array}{cc}
      Q_0^r[B_0;B_0] & Q_0^r[B_0;H'_0] \\
         0  & Q_0^r[H'_0;H'_0] \\
          \end{array} \right)

   \end{array},  \]
they agree on the diagonal blocks and so the diagonal blocks of
$T_\gamma Q_0^r$ and those of $T_\gamma \overline{Q}_0^r$ are the
same. Hence these finite rank operators have the same trace. \\

For $d \geq q \geq 1$ the complex $(\mathcal{C})$ is exact at the
$q^{th}$ place, hence $C_q = B_q \oplus B'_{q-1}$ where $B'_{q-1}
\overset{\partial}{\underset{\sim}{\longrightarrow}}
B_{q-1}$ and $\oplus$ is $T_\gamma$-stable.\\

Since for $d \geq q \geq 1$, $Q_q^r$ and $\overline{Q}_q^r$ differ
also on the bottom (right) diagonal block:

\[ \begin{array}{cc}

 Q_q^r = \left( \begin{array}{cc}
      Q_q^r[B_q;B_q] & Q_q^r[B_q;B'_{q-1}] \\
      Q_q^r[B'_{q-1};B_q] & Q_q^r[B'_{q-1};B'_{q-1}] \\
          \end{array} \right)
 &
 \overline{Q}_q^r = \left( \begin{array}{cc}
      Q_q^r[B_q;B_q] & Q_q^r[B_q;B'_{q-1}] \\
         0  & \widetilde{Q}^r_{q-1} \\
          \end{array} \right)

   \end{array}  \]
we cannot conclude yet that $T_\gamma Q_q^r$ and $T_\gamma
\overline{Q}_q^r$ have the same diagonal blocks with respect to
$C_q =B_q \oplus B'_{q-1}$. We now use the (refined) nice $\oplus$
decomposition of $C_q$. Recall that
$(B_{q-1}^{in(r)})'=\alpha(B_{q-1}^{in(r)})$ and that
$(B_{q-1}^{out(r)})'=\alpha(B_{q-1}^{out(r)})$, where $\alpha$ is
the section of $\partial$ of Lemma ~\ref{key lemma}.

With respect to the decomposition $ C_q = B_q \oplus
(B_{q-1}^{in(r)})' \oplus (B_{q-1}^{out(r)})' $ we have:
\[ T_\gamma = \left( \begin{array}{ccc}
      * & 0 & 0\\
      0 & * & 0\\
      0 & 0 & *
          \end{array} \right).            \]

Since $Q_q^r$ and $\overline{Q}_q^r$ agree on $B_q$ we show that
they have the same diagonal blocks with respect to
$(B_{q-1}^{in(r)})'\oplus (B_{q-1}^{out(r)})'$.

Since $\overline{Q}_q^r[B'_{q-1};B'_{q-1}] =
\widetilde{Q}_{q-1}^r= Q_{q-1}^r[B_{q-1};B_{q-1}]$ we look at
$Q_{q-1}^r$ on $B_{q-1} = B_{q-1}^{in(r)} \oplus
B_{q-1}^{out(r)}$.\\

Note that:
\begin{eqnarray*}
B_{q-1}^{out(r)} \cap C_{q-1}^{in(r)}& = &(B_{q-1}^{out(r)} \cap
B_{q-1}) \cap C_{q-1}^{in(r)}\\& = &B_{q-1}^{out(r)} \cap (B_{q-1}
\cap C_{q-1}^{in(r)})\\ &=& B_{q-1}^{out(r)} \cap B_{q-1}^{in(r)}
= \{0\}.
\end{eqnarray*}

Now for $b \in B_{q-1}^{out(r)}$, have
$Q_{q-1}^r[B_{q-1}^{out(r)}; B_{q-1}^{out(r)}] (b) \in
B_{q-1}^{out(r)} \cap C_{q-1}^{in(r)} = \{0\}$, hence $Q_{q-1}^r$
acts as the zero operator on $B_{q-1}^{out(r)}$.

On $B_{q-1}^{in(r)}$ we have $Q_{q-1}^r$ acting as the identity
operator since $B_{q-1}^{in(r)} \subset C_{q-1}^{in(r)}$. So
$\overline{Q}_q^r$ with respect to $C_q = B'_{q} \oplus
(B_{q-1}^{in(r)})' \oplus (B_{q-1}^{out(r)})'$ has the form:

\[ \overline{Q}_q^r = \left( \begin{array}{ccc}
      Q_q^r & * & *\\
      0 & 1 & *\\
      0 & 0 & 0
          \end{array} \right).            \]
We now find the form of $Q_{q}^r$ with respect to $C_q = B'_{q}
\oplus (B_{q-1}^{in(r)})' \oplus (B_{q-1}^{out(r)})'$.\\

Let $b \in (B_{q-1}^{out(r)})'$, then $\partial b \in
B_{q-1}^{out(r)}$. $Q_{q}^r b = b_q + b^{in} + b^{out} \in
C_q^{in(r)}$ so $\partial Q_{q}^r b = 0 + \partial b^{in} +
\partial b^{out} \in B_{q-1}^{in(r)}$. Hence $\partial b^{out}
= 0$ and so $b^{out}= 0$ which means that
$Q_{q}^r[(B_{q-1}^{out(r)})';(B_{q-1}^{out(r)})'] = 0$. Now,
$(B_{q-1}^{in(r)})'\subset C_{q-1}^{in(r)}$ so $Q_{q}^r = 1$ on
$(B_{q-1}^{in(r)})'$. Hence we obtain that $Q_{q}^r$ has the form:

\[ Q_q^r = \left( \begin{array}{ccc}
      Q_q^r & * & *\\
      * & 1 & *\\
      * & 0 & 0
          \end{array} \right)           \]

Thus $Q_q^r$ and $\overline{Q}_q^r$ have the same diagonal blocks
and so $trace(T_\gamma Q_q^r) = trace(T_\gamma \overline{Q}_q^r)$.
\end{proof}

\section{\label{SS:SL_2}Recovering Schneider-Stuhler's result for semisimple
rank $1$ groups}

The following lemma is part of the proof of lemma 12 in [S-S
III.4]. It holds for any connected reductive group $G$ with
building $X$. For the sake of completeness we repeat the proof
here.
\begin{lem}\label{ellopennbd} Let $\gamma \in G^{ell}$.
There exists an open subgroup $U \subset G$ such that
$$ X^\gamma =X^{\gamma '} \; \mbox{ for all } \; \gamma ' \in \gamma U .$$
\end{lem}
\begin{proof}
Since $\gamma$ is elliptic, its fixed point set is non-empty and
compact. Fix a point $o \in X^\gamma$ and choose $r \gneq 0$ large
enough so that $X^\gamma$ is contained in the interior of the
closed ball $B(o,r)$. By construction $B(o,r)$ has no
$\gamma$-fixed points on its boundary. Let $U$ denote the open
subgroup: $U:=\{g \in G \; | \; gx=x \; \mbox{ for all } \; x\in
B(o,r)\}$. For $\gamma' \in \gamma U$ we show that $X^\gamma =
X^{\gamma '}$.\\
It is clear that the actions of $\gamma$ and $\gamma '$ agree on
all the points of the closed ball $B(o,r)$, so in particular they
have the same fixed points inside $B(o,r)$: $X^\gamma \cap B(o,r)
= X^{\gamma '}\cap B(o,r)$. Suppose $\gamma '$ has a fixed point
$x$ outside the ball. Then the whole geodesic $geod(o,x)$ must be
fixed by $\gamma '$ and so there is a $\gamma '$-fixed point
$x_0:=\partial B(o,r)\cap geod(o,x)$ on the boundary of the ball.
But any such $\gamma '$-fixed point is also a $\gamma$-fixed
point, contradiction $\gamma$ having no fixed points on the
boundary of $B(o,r)$.
\end{proof}
Let $(\pi,V)$ be an admissible representation of $G$ with
character function $\Theta_\pi(g)$, defined on the regular
elements in $G$. Let $\gamma$ be a regular semisimple elliptic
element: $\gamma \in G^{ell}$, with $o \in X^\gamma$, and fix $e
\geq e_0$. We recover the character formula ~(\ref{chi}) of
Schneider--Stuhler.\\
\begin{claim}\label{Kgroup}
We can fix an open compact subgroup $K$ of $G$ with the following
properties:
\begin{enumerate}
\item{} The character is locally constant on the neighborhood $\gamma K$ of
$\gamma$, i.e. $\Theta_\pi(\gamma)= \Theta_\pi(\gamma K)$
\item{} All elements $\gamma k$, $k \in K$, have the same fixed point set:
$X^\gamma = X^{\gamma k}$
\item{} $K \subset \underset{x\in X^\gamma}\bigcap U_x^{(e)}$
\item{} $\gamma$ normalizes $K$.
\end{enumerate}
\end{claim}
\begin{proof}
To see that such a group $K$ exists, consider the following. If
$U$ denotes the open subgroup of the lemma above, and $U_x^{(e)}$
the usual open (compact) subgroups attached to a point $x\in X$,
then $U \underset{x\in X^\gamma}\bigcap U_x^{(e)}$ is an open
neighborhood of the identity in $G$. Since $\gamma$ fixes $o$ it
is contained in the group $P_o^{\dagger}$, and so it normalizes
all the open compact groups $U_o^{(e)}$, (property $(U2)$ above).
Thus if we let $K:=U_o^{(e)}$ and choose $e$ large enough, we can
make sure that $\Theta_\pi(\gamma K)= \Theta_\pi(\gamma)$ and $K
\subset U \underset{x\in X^\gamma}\bigcap U_x^{(e)}$. Such $K$
satisfies all the properties above.
\end{proof}
Set $f:=\frac {1_{\gamma K}}{vol(K)} \in C^{\infty}_c(G)$ to be
the characteristic function of the set $\gamma K$ normalized by
its volume. We have:
\begin{eqnarray*}
trace(\pi(f);V) & = & \underset{K}\int f(\gamma k)trace(\pi(\gamma k);V)\,dk\\
                & = & \underset{K}\int \frac{1}{vol(K)}\Theta_\pi(\gamma k) \, dk\\
                & = & \Theta_\pi(\gamma) \underset{K}\int \frac{1}{vol(K)} \, dk\\
                & = & \Theta_\pi(\gamma)
\end{eqnarray*}

Choose $r=r(f,\gamma)$ large enough so that is satisfies the
requirement of Theorem~\ref{thmfunction} and so that $X^\gamma \subset B(o,r)$.\\

For the rest of this section $G$ will denote a group of semisimple
rank $1$. Recall that for such groups we have shown (Corollary
~\ref{directsumexist}) that nice $\oplus$ decompositions exist
uniformly in $e$ (with respect to $\gamma \in G^{cpt}$ such that
$o\in X^\gamma$).\\

All elements $\gamma k \in \gamma K$ share the same fixed point
$o$. We claim that it is possible to choose a nice decomposition
$(\oplus)$ which is common to all such elements: the issue being
that we want the decomposition to be $T_{\gamma k}$-equivariant
simultaneously for all $ k\in K$.\\

In the construction of a nice decomposition ($\oplus$) we averaged
over the group $K_\gamma$ so as to make sure the construction was
$\gamma$-equivariant. Now, we will average over $K_\gamma$ and
then average again over $K$. Property $(4)$ of the group $K$
implies that the construction is both $\gamma$-equivariant and
$K$-equivariant. We demonstrate this process with the section
$\alpha$. Start with $\alpha$ and average over $K_\gamma$ to
obtain:

$$\alpha':=\int\limits_{K_\gamma} T_k^{-1} \circ \alpha \circ T_k$$

which is $\gamma$-equivariant. Now average $\alpha'$ over $K$ the
obtain:

$$\alpha'':=\int\limits_{K} T_k^{-1} \circ \alpha' \circ T_k$$

By construction, $\alpha''$ is $K$-equivariant. We check that it
is still $\gamma$-equivariant:
\begin{eqnarray*}
T_\gamma^{-1} \circ \alpha'' \circ T_\gamma & = & \int\limits_{K}
T_\gamma^{-1}
\circ T_k^{-1} \circ \alpha' \circ T_k \circ T_\gamma\\
  & = & \int\limits_{K} (T_\gamma^{-1}
\circ T_k^{-1} \circ T_\gamma) \circ (T_\gamma^{-1} \circ \alpha'
\circ T_\gamma)
\circ (T_\gamma^{-1} \circ T_k \circ T_\gamma)\\
& =& \int\limits_{K} T_{k'}^{-1} \circ \alpha' \circ T_{k'} \; =
\; \alpha''
\end{eqnarray*}

Thus the nice decomposition works uniformly on $\gamma K$ and so
by Proposition~\ref{ridofbar}:
\[trace(T_{\gamma k} Q_q^r;C_q) = trace(T_{\gamma k} \overline{Q}_q^r;C_q)
\quad d \geq q \geq 0 \qquad \mbox{for all} \; k \in K.
\]
For $\gamma$ regular elliptic $trace(T_{\gamma}; C_q)$ is
essentially counting stable $q$-facets (with multiplicity) and
hence we see that:
$$trace(T_{\gamma}; C_q)= \underset{\gamma-stable \atop F\in X_q }
\sum {(-1)}^{q-dim F(\gamma)}trace(\gamma; V^{U_F^{(e)}}).$$ Here
the quantity $q-dim F(\gamma)$ is $\pm 1$, depending on whether
$\gamma$ preserves or reverses the orientation of $F$ (see
~\cite[p.45 and p.51]{SS1} for details).

In the notation of ~\cite{SS2} we have:
\begin{eqnarray*}
\sum_{q=0}^{d} (-1)^q trace(T_{\gamma}; C_q) & = & \sum_{q=0}^{d}
\underset{\gamma-stable \atop F\in X_q } \sum {(-1)}^{dim
F(\gamma)}trace(\gamma;
V^{U_F^{(e)}})\\
& = & \sum_{q=0}^{d} \underset{F(\gamma) \in {(X^\gamma)}_q}\sum
{(-1)}^{q}trace(\gamma; V^{U_F^{(e)}})
\end{eqnarray*}
Among other things this says that the left-hand-side equation is
constant on $\gamma K$ (since the right-hand-side is the same for
all $\gamma k$, $k \in K$, by the choice of $K$).\\

Using this equality and Theorem~\ref{thmfunction} (with $f:=\frac
{1_{\gamma K}}{vol(K)}$ and the $r$ above) we obtain:
\begin{eqnarray*}
\sum_{q=0}^{d} (-1)^q trace(T_{\gamma}; C_q)\underset{K}\int \, dk
& = & \underset{K}\int \left\{ \sum_{q=0}^{d} (-1)^q
trace(T_{\gamma k}; C_q)\right\} \, dk \\ & = & \underset{K}\int
\left\{ \sum_{q=0}^{d} (-1)^q trace(T_{\gamma k}{Q}_q^r;
C_q)\right\} \, dk \\& = & \underset{K}\int
\left\{ \sum_{q=0}^{d} (-1)^q trace(T_{\gamma k}\overline{Q}_q^r; C_q)\right\} \, dk\\
& = &
      \sum_{q=0}^{d}  (-1)^q trace(T_{f}\overline{Q}_q^r ; C_q)\\
      & = &  trace(\pi(f) ; V) = \Theta_\pi(\gamma).
\end{eqnarray*}

Putting all of the above together we recover the Schneider-Stuhler
formula ~(\ref{chi}), for groups of semisimple rank $1$:
\begin{thm} For $G$ a semisimple rank $1$ group , $V$ a finitely generated
admissible representation of $G$, $\gamma \in G$ regular elliptic
and $e \geq e_0(V)$:
\[ \Theta_\pi(\gamma) = \sum_{q=0}^{d}  \underset{F(\gamma) \in
{(X^\gamma)}_q}\sum (-1)^q trace(\gamma; V^{U_F^{(e)}})\]
\end{thm}

{\bf Note:} unlike the Schneider-Stuhler proof of this theorem,
our proof does not rely on Kazhdan's density theorem.
\section{A character formula for compact elements}
In this section we prove the main results of this paper: Theorem
~\ref{main} and Corollary ~\ref{maincor}. The key fact used in
proving the main result, which will be developed here, is that a
regular semisimple compact element $\gamma$ has an open
neighborhood in $G$ such that, as we vary $\gamma'$ inside such
neighborhood, the following expression stays constant:
\begin{eqnarray}\label{heartsuit} \sum_{q=0}^{d} \underset{F(\gamma')\in(X^{\gamma '} \cap X^r)_q}
\sum {(-1)}^q trace(\gamma ';V^{U_{F}^{(e)}})
\end{eqnarray}

\subsection{The periodic nature of $X^\gamma$}
\begin{figure}[h]
\psfragscanon \centering \psfrag{o}{$o$}
\psfrag{Finite}{$\mbox{finite information}$}
\epsfig{file=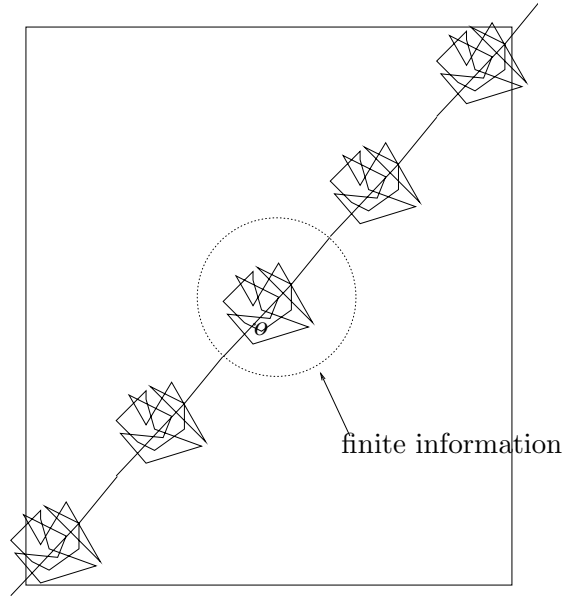,height=8cm}\caption{\label{periodic}Schematic
description of the periodic nature of $X^\gamma$}
\end{figure}

When $\gamma$ is elliptic, we know that its fixed point set,
$X^{\gamma}$, is finite (as a simplicial complex). When $\gamma$
is compact, the fixed point set, $X^{\gamma}$, can be an infinite
simplicial complex. Yet, as we now show, this set has periodic
nature. [See Figures ~\ref{periodic} and ~\ref{tree1}]. Recall
that the centralizer
of $\gamma$, $C_G(\gamma)$, acts on $X^\gamma$.\\

Let $K=K_o$ denote the stabilizer of the vertex $o\in X^\gamma$.
$K$ is an open subgroup of $G$, which is compact modulo the center
of $G$. In fact it is a maximal such subgroup. Denote its
characteristic function by $ch_K$.\\

We now assume for simplicity that all of the vertices of the
building $X$ have the same type. This means that we can identify
the set of cosets $G/K$ with the vertices of the building, $X_0$,
via $gK \rightarrow go$. (Without this assumption the vertices can
have a finite number of types and a little more care is needed in
keeping track of the different types.)\\

Fix a Haar measures $dg$ on $G$ and $dt$ on $T$. Let
$\frac{dg}{dt}$ be the invariant measure on $T\backslash G$ with
respect to $dg$ and $dt$. The orbital integral
$$\mathcal{O}_{\gamma}(ch_K):=\underset{T \backslash G}\int ch_K(g^{-1}
\gamma g)\, \frac{dg}{dt}$$ where $T:=C_G(\gamma)^{\circ}$ is
known to
converge, say to the (finite) number $N$.\\

\begin{rem} An orbital integral such as $\underset{G}\int ch_K(g^{-1}
\gamma g)\, dg$ does not converge in general. The issue is that
the split part of $T$, might not be compact and may make this
integral diverge. In our case we integrate over $T\backslash G$ so
we are avoiding this type of problem. Also note that since $T$ is
generally of smaller dimension than $G$, we need two normalization
factors, one for $dg$ and one for $dt$.
\end{rem}

Since a vertex $gK \in X$ is fixed by $\gamma$ ($\gamma gK=gK$) if
and only if $g^{-1} \gamma g \in K$, we have the following:

\begin{eqnarray*}
N = \underset{T \backslash G}\int ch_K(g^{-1} \gamma g)\,
\frac{dg}{dt} & = &
    \underset{g\in T \backslash G / K}\sum \frac {vol_G(K)}{vol_T(T \cap K)} ch_K(g^{-1} \gamma g)\\
& = & \underset{{g\in T \backslash G / K} \atop {\gamma gK=gK}}
\sum \frac {vol_G(K)}{vol_T(T \cap K)}\\
& = & \frac {vol_G(K)}{vol_T(T \cap K)} (\mbox{number of T-orbits
on }(G/K)^{\gamma}).
\end{eqnarray*}

Thus the set of $\gamma$-fixed vertices,
$X_0^{\gamma}=(G/K)^{\gamma}$, is the union of a finite number of
$T$-orbits. Similarly, if we allow vertices of different types and
consider also higher dimensional facets, we reach the same
conclusion: that the fixed point set $X^\gamma$ is the union of a
finite number of $T$-orbits (Figure ~\ref{periodic} shows one such
$T$-orbit). That is, there is a \textit{fundamental domain} $D$ in
$X$, which is a \textit{finite} union of facets, such that any
facet in $X^\gamma$ is of the form $tF$, for some facet $F\in D$
and some $t\in T$ (see Figure ~\ref{tree2}). This is what we mean
by saying that `the set $X^\gamma$ has periodic nature'.

\begin{rem} In the case that $\gamma$ is elliptic, so that $T$ is
compact, each $T$-orbit on $X^{\gamma}_0$ contains finitely many
points. But in the case that $\gamma$ is compact but not elliptic,
each $T$-orbit is infinite. (Essentially the maximal $k$-split
torus contained in $T$, acts as a translation.)
\end{rem}

\subsection{An open neighborhood of $\gamma$ whose elements share the same
number of fixed points} In the elliptic case, Lemma
~\ref{ellopennbd} guaranteed the existence of an open neighborhood
$\gamma U$ of $\gamma$ such that $ X^\gamma=X^{\gamma'} \; \mbox{
for all } \; \gamma' \in \gamma U.$ We used the existence of such
a neighborhood in showing that the elements $\gamma '=\gamma u$,
$u\in U$, all had the same (number of) fixed facets inside any
ball $B(o,r)$ of radius $r \geq 0$ about $o \in X^\gamma$. In
general, for a compact element it is not possible to find an open
neighborhood $\gamma U$ whose elements share the same fixed point
set. In this section we show that it is possible to find an open
neighborhood $\gamma U$ whose elements all share the same
\emph{number} of fixed points inside any such ball $B(o,r)$.\\

To understand the behavior of elements near $\gamma$ it is
convenient to assume $\gamma$ to be regular semisimple (in fact
until now it was not necessary to assume $\gamma$ to be
semisimple).\\

For any maximal torus $T$, let $T'$ denote the set of regular
elements of $T$. Recall Harish-Chandra's \emph{submersion
principle}, ~\cite[Lemma 20, p.55]{Harish-Chandra}, which says
that the following map is submersive:

\begin{eqnarray*}
\psi \; : \; G/T \times T' & \longrightarrow &  G\\
   (gT\;,\; t) & \rightarrow & g t g^{-1}
\end{eqnarray*}

This implies that the image of this map $$\psi (G/T\times T')=\{g
t g^{-1}\;|\; g\in G, \; t \in T'\}=\mathcal{O}_G(T'),$$ is an
open subset in $G$. Given $\gamma\in G^{reg}$, let
$T=C_G(\gamma)^\circ$. In this case $\psi (G/T\times T')$ is an
open neighborhood of $\gamma$ in $G^{reg}$. So in a small
neighborhood of $\gamma$, the torus $T$ and the conjugacy class
$\mathcal{O}_G(\gamma):=\{g \gamma g^{-1} \;|\;
g\in G\}$ are transversal and $\mathcal{O}_G(T')$ fills out an open neighborhood of $\gamma$.\\

\begin{figure}[!h]
\psfragscanon \centering \psfrag{h}{$\gamma$} \psfrag{T}{$T$}
\psfrag{T'}{$T^{\circ}$}  \psfrag{O(h)}{$\mathcal{O}_G(\gamma)$}
\epsfig{file=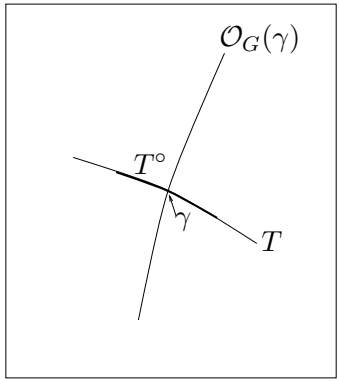,height=5cm}\caption{\label{transversal}
Local picture around $\gamma$}
\end{figure}

\begin{figure}[h]
\psfragscanon \centering \psfrag{h}{$\gamma$} \psfrag{T}{$T$}
\psfrag{O(h)}{$\mathcal{O}_G(\gamma)$} \psfrag{T'}{$T^{\circ}$}
\psfrag{hK}{$\gamma K$}
\psfrag{O(T,U)}{$\mathcal{O}_G(T^{\circ})$}
\epsfig{file=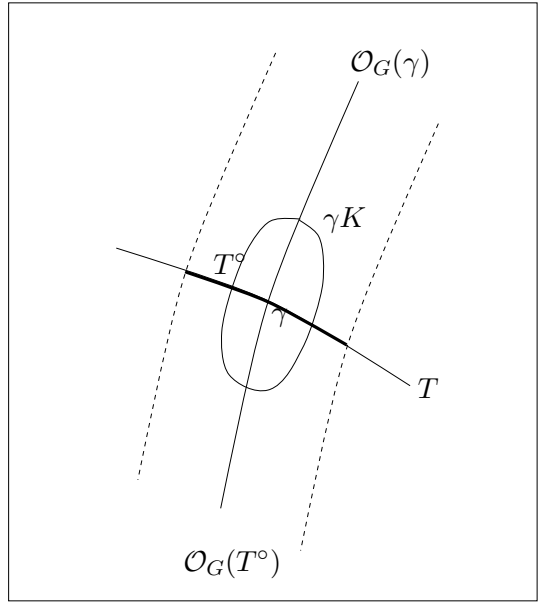,height=8cm}\caption{\label{tubular}Tubular
neighborhood of $\gamma$}
\end{figure}

This means that to understand an open neighborhood of (the regular
semisimple) $\gamma$ in $G$, it is enough to understand an open
neighborhood of $\gamma$ in inside $T:=C_G(\gamma)^{\circ}$ and an
open neighborhood of $\gamma$ inside the conjugacy class
$\mathcal{O}_G(\gamma)$ [See Figures ~\ref{transversal} and
~\ref{tubular}].\\

We first study the behavior of elements in an open neighborhood
of $\gamma$ in $T$.\\

Let $D$ be a fundamental domain for the action of $T$ on
$X^\gamma$. Fix a point $o\in D$ and choose $r \gneq 0$ large
enough so that $D$ is contained in the interior $B^\circ$ of the
closed ball $B=B(o,r)$. Let $U$ denote the open subgroup: $U:=\{g
\in G \; | \; gx=x \; \mbox{ for all } \; x\in B(o,r)\}$, the
pointwise stabilizer of $B(o,r)$. We will show that any element
$\gamma ' \in T^{\circ}:=T \cap \gamma U$ has the same fixed set
as $\gamma$: $X^\gamma = X^{\gamma '}$.\\

\begin{figure}[h]
\psfragscanon \centering \psfrag{Xh}{$X^{\gamma}$}
\psfrag{Xg}{$X^{g\gamma' g^{-1}}$} \psfrag{o}{$o$}
\psfrag{O(h)}{$\mathcal{O}_G(\gamma)$}
\epsfig{file=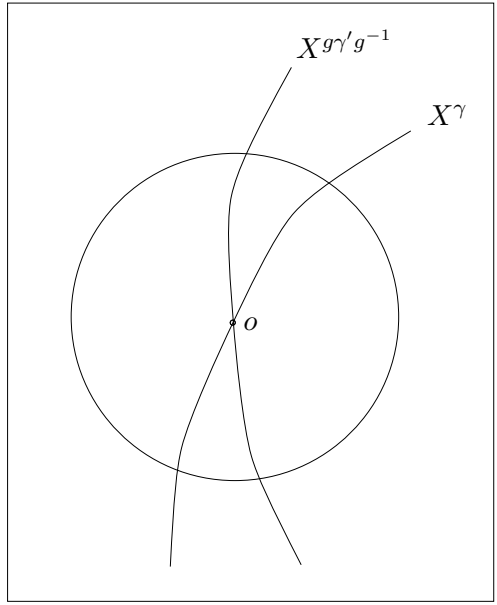,height=8cm}\caption{\label{rotated} The
fixed point sets $X^{\gamma}$ and $X^{g\gamma' g^{-1}}$}
\end{figure}

Any facet $F \in X^\gamma$ belongs to a $T$-orbit: say $F\in TF_0$
for some $F_0 \in D$; so there exists a $t\in T$ such that $tF =
F_0$. Since any element $\gamma '\in T^{\circ}$ fixes $F_0$ and
commutes with all elements of $T$, we have:
$$ \gamma 'F = \gamma't^{-1}tF = t^{-1}\gamma'(tF) = t^{-1}\gamma'F_0
= t^{-1}F_0 = F.$$ We see that $\gamma'$ fixes all the facets in
$X^\gamma$; that is for all $\gamma' \in T^{\circ}$: $X^\gamma
\subset X^{\gamma'}$.
\begin{figure}[h]
\psfragscanon \centering  \psfrag{o}{$o$}
\psfrag{T^o}{$\mathcal{T}^o$} \psfrag{D}{$D$} \psfrag{B^o}{$B^o$}
\epsfig{file=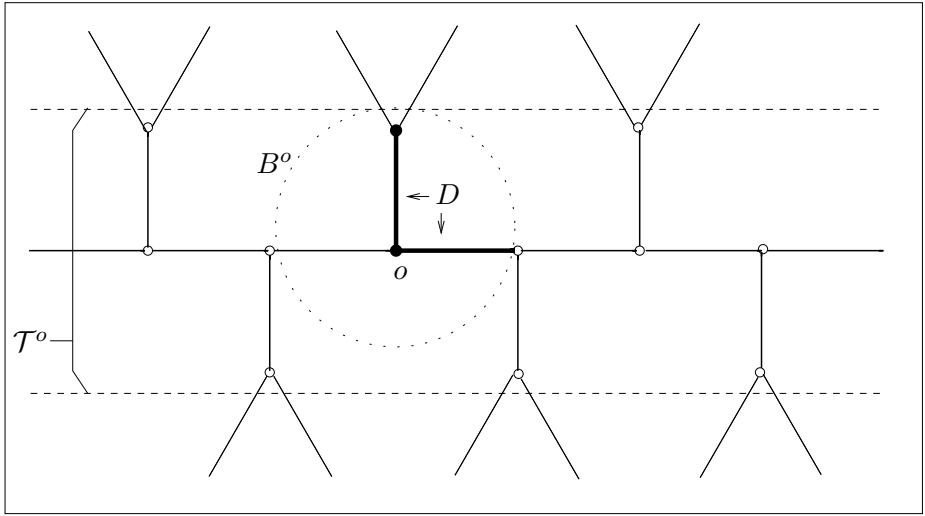,height=8cm}\caption{\label{tree2}A
fundamental domain (in bold) $D$ in the fixed point set of the
compact non-elliptic element $\gamma\in GL_2$ of Example
~\ref{exnce}. }
\end{figure}
\begin{ex}
We show that a compact element $\gamma \in G=GL_2$ has an open
neighborhood $T^\circ$ in $T=C_G(\gamma)^\circ$, such that all
elements in $T^\circ$ have the same fix point set.\\

Let $G=GL_2$, and take a compact element $\gamma= \left(
\begin{array}{cc}
               1   &   0 \\
               0   &   u
           \end{array} \right)$,
where $u$ has the form: $\;u= 1+\alpha_r\varpi^r+\cdots  \quad$
$(0 \neq \alpha_r \in \overline{k}).$ Recall (Example
~\ref{exnce}) that the points of $X$ fixed by $\gamma$ are
precisely those at a distance less than or equal to $r$ from
the basic apartment $A$.\\

Consider the following open neighborhood of $u$: $u+\varpi^{r+1}
O= 1+\alpha_r\varpi^r+\varpi^{r+1} O$. All elements in this
neighborhood have the same valuation. Elements in  the
corresponding open neighborhood of $\gamma$ in $T$:
$$T^\circ = \{ \left( \begin{array}{cc}
               1   &   0 \\
               0   &   u'
           \end{array} \right) | u' \in u+\varpi^{r+1} O \}$$
all have the same fixed point set.
\end{ex}
We show this is the case in general.
\begin{lem}\label{fixedconstant}
All elements $\gamma'\in T^{\circ}:=T \cap \gamma U$ have the same
fixed point set: $X^{\gamma'}=X^\gamma$.
\end{lem}
\begin{proof}
We imitate the proof of Lemma ~\ref{ellopennbd} which worked in
the elliptic case. Using the notation above, let the open ball
$B^{\circ}$ be the interior of the closed ball $B$. Let
$\mathcal{T}^{\circ}:=TB^{\circ}$ be the $T$-orbit $B^{\circ}$.
Since the subsets $tB^{\circ}$ are open and since $TB^{\circ}=
\bigcup_{t\in T} tB^{\circ}$, we see that $\mathcal{T^{\circ}}$ is
an open subset of $X$. We will call $\mathcal{T^{\circ}}$, the
\textit{open tube} around $X^\gamma$, and its closure,
$\mathcal{T}:=\overline{\mathcal{T^{\circ}}}$, the \textit{closed
tube} around $X^\gamma$. Figure ~\ref{tree2} shows these sets for
the compact non-elliptic element $\gamma$ (with $r=1$) of Example
~\ref{exnce}. Note that $B^{\circ} \supset D$ and hence
that $\mathcal{T}^{\circ}= TB^{\circ} \supset TD = X^\gamma$.\\
It is easy to see that $\mathcal{T}=TB$.\\
Now let $\gamma' \in \gamma U$. Since the actions of $\gamma$ and
$\gamma'$ agree on the closed ball $B$, they also agree on the
closed tube $\mathcal{T}=TB$:
$$\gamma t b = t \gamma b = t \gamma' b = \gamma' t b \qquad
t\in T \quad b\in B.$$ Since $X^\gamma \subset
\mathcal{T^{\circ}}$, $\gamma$ has no fixed points on the
boundary, $\partial \mathcal{T}$, of $\mathcal{T}$. Hence also
$\gamma'$ has no fixed points on
$\partial \mathcal{T}$.\\
Now suppose $\gamma'$ has a fixed point $x$ outside the closed
tube $\mathcal{T}$. Then the whole geodesic $geod(o,x)$ must be
fixed by $\gamma '$ and so there is a $\gamma '$-fixed point
$x_0:=\mathcal{T}\cap geod(o,x)$ on the boundary of the closed
tube $\mathcal{T}$, contradicting $\gamma'$ having no fixed points
on the boundary of $\mathcal{T}$. Thus $X^{\gamma'}\subset
\mathcal{T}$, and hence $X^{\gamma'} = X^{\gamma}$.
\end{proof}
\newpage
Elements in the tubular neighborhood $\mathcal{O}_G(T^{\circ})$ of
$\gamma$ all look like $g \gamma' g^{-1}$, $\gamma' \in T^{\circ},
g\in G$. [See Figure ~\ref{tubular}]. For such elements and $x\in
X^{\gamma}$, we have: $g\gamma' g^{-1}(gx) = g\gamma 'x = gx$.
That is $X^{g\gamma' g^{-1}} = gX^{\gamma'} = g X^{\gamma}$. If we
choose $g$ in a small enough neighborhood $U^{\circ}$ of the
identity in $G$ (e.g. $g\in U_{o}^{(e)}$) so that it fixes $o$,
then $X^{g\gamma'g^{-1}}$ looks like a rotated version of
$X^\gamma$. [See Figure ~\ref{rotated}]. Even though the fixed
point set $X^{g\gamma' g^{-1}}$ is not the same as $X^\gamma$,
inside any ball $B(o,r)$, $X^{g\gamma' g^{-1}}$ and $X^\gamma$
have the same number of fixed facets. Since the construction of
$X^r$ is $g$-equivariant, the same is true inside any truncated
building
$X^r$.\\
That is, since elements in a small enough open neighborhood of
$\gamma$ are all of the form $g\gamma 'g^{-1},\; \gamma' \in
T^{\circ},\, g\in U^{\circ}$, we conclude that there exists an
open neighborhood of $\gamma$ in $G$ whose elements all have the
same number of fixed points inside any ball $B(o,r)$, hence inside
any truncated building $X^r$ about $o$.
\subsection{Traces of elements in $T$ close to $\gamma$}
We now study the traces (on the fibers over their fixed points) of
elements $\gamma' \in T$ close to $\gamma$. We already saw that
for elements $\gamma' \in T^{\circ}$: $X^{\gamma} = X^{\gamma'}$.
Now we show that for $\gamma'\in T^{\circ}$ and $x\in
X^{\gamma'}=X^\gamma$, $trace(\gamma',V^{U_x^{(e)}})$ is constant
on the $T$-orbit of $x$.
\begin{lem} Given $\gamma' \in T^{\circ}$, $e\geq 0$, and
$x\in X^{\gamma '}=X^\gamma$. For any $s\in T$ we have:
$$ trace(\gamma',V^{U_x^{(e)}}) = trace(\gamma',V^{U_{sx}^{(e)}}) $$
\end{lem}
\begin{proof}
Since $U_{sx}^{(e)} = s U_x^{(e)} s^{-1}$ we have the isomorphism:
\begin{eqnarray*}
\phi \; : \; V^{U^{(e)}_x} & \longrightarrow & V^{U^{(e)}_{sx}}\\
   v & \mapsto & sv
\end{eqnarray*}
and using the fact that $T$ is abelian we get that for any
$\gamma' \in T$ the following diagram commutes:
\begin{eqnarray*}
     v       & \overset{\phi}\longrightarrow &     sv\\
\gamma' \downarrow &                          & \downarrow \gamma'\\
    \gamma' v       & \overset{\phi}\longrightarrow & \gamma' sv=s\gamma' v
\end{eqnarray*}
Hence the $trace$ of $\gamma'$ on $V^{U_x^{(e)}}$ is the same as
its $trace$ on $V^{U_{sx}^{(e)}}$.
\end{proof}

\begin{cor}
For $\gamma$ a (regular semisimple) compact element there is an
open neighborhood $T^{\circ \circ}\subset T^{\circ}$ of $\gamma$
in $T=C_G(\gamma)^{\circ}$ (which depends on $\gamma$ and $e$)
such that all elements in this neighborhood have the same fixed
point set, $X^{\gamma}$, and they all act the same way on the
fibers above each fixed point $x\in X^{\gamma}$:
$$ trace(\gamma',V^{U_x^{(e)}}) = trace(\gamma,V^{U_{x}^{(e)}}) \qquad
\mbox{for all } \gamma' \in T^{\circ \circ}$$
\end{cor}
\begin{proof}
By the last lemma, for $\gamma' \in T^{\circ}$, the traces
$trace(\gamma',V^{U_x^{(e)}})$ and\\
$trace(\gamma,V^{U_x^{(e)}})$, are constant on the $T$-orbit of
$x$. Hence it is enough to show this equality for representatives
of the orbits. Let $x\in D$ and let $\bigcap_{x\in D} U_{x}^{(e)}$
be an open neighborhood of the identity in $G$ of elements which
act trivially on all the fibers $V^{U_x^{(e)}}$, $x\in D$. Then
all elements $\gamma'\in \gamma ( \bigcap_{F\in D} U_{F}^{(e)})$
act on fibers the same way as $\gamma$, so that for $\gamma' \in
T^{\circ \circ}:=\gamma U \cap \gamma ( \bigcap_{F\in D}
U_{F}^{(e)})$, we have: $trace(\gamma' ,V^{U_x^{(e)}}) =
trace(\gamma,V^{U_{x}^{(e)}})$.
\end{proof}
\begin{cor}
The expression: $$ \sum_{q=0}^{d}
\underset{F(\gamma')\in(X^{\gamma '} \cap X^r)_q} \sum {(-1)}^q
trace(\gamma ';V^{U_{F}^{(e)}})$$ is constant for all $\gamma'\in
T^{\circ \circ}$.
\end{cor}
\begin{proof}
By the last corollary $trace(\gamma ';V^{U_{F}^{(e)}})$ is
constant for all $\gamma' \in T^{\circ \circ}$. Since $X^{\gamma
'}$ is constant for all $\gamma' \in T^{\circ \circ}$ by Lemma
~\ref{fixedconstant}, the expression $(X^{\gamma '} \cap X^r)_q$
is also constant for such $\gamma'$ and so the summation is over
the same facets $F$ as $\gamma'$ varies. Hence the whole
expression above is constant for $\gamma'\in T^{\circ \circ}$.
\end{proof}
\subsection{Traces of elements in $\mathcal{O}_G(T^{\circ\circ})$ close to $\gamma$}
In this subsection we show that the constancy of the expression
above (same as expression ~(\ref{heartsuit})) on the open
neighborhood $T^{\circ\circ}$ of $\gamma$ inside $T$ extends to a
open neighborhood of $\gamma$ inside $G$.

Recall that for elements $g$ in the small neighborhood of the
identity and $\gamma' \in T^{\circ}$, we have: $X^{g\gamma'
g^{-1}}= gX^{\gamma'}=gX^{\gamma}$.

Now, let $x$ be a fixed point of $\gamma' \in T^{\circ\circ}$ and
let $gx$ the corresponding fixed point of $g\gamma' g^{-1}$
(recall Figure ~\ref{rotated}). Note that if $g$ fixes the origin
$o$, then we have that $x\in X^r$ if and only if $gx \in X^r$.

\begin{lem}
With the notation above, for any $\gamma' \in T^{\circ}$, any $g
\in U^\circ$, and any $e \geq 0$ we have:
$$ trace(\gamma',V^{U_x^{(e)}}) = trace(g\gamma' g^{-1},V^{U_{gx}^{(e)}}) $$
\end{lem}
\begin{proof}
This proof is almost identical to the proof of the previous Lemma.
Since $U_{gx}^{(e)} = g U_x^{(e)} g^{-1}$ we have the isomorphism:
\begin{eqnarray*}
\phi \; : \; V^{U^{(e)}_x} & \longrightarrow & V^{U^{(e)}_{gx}}\\
   v & \mapsto & gv
\end{eqnarray*}
for any $\gamma' \in T^{\circ}$ the following commutative diagram:
\begin{eqnarray*}
     v       & \overset{\phi}\longrightarrow &     gv\\
\gamma' \downarrow &                          & \downarrow g\gamma' g^{-1}\\
    \gamma' v       & \overset{\phi}\longrightarrow & g\gamma' v
\end{eqnarray*}
Hence the $trace$ of $\gamma'$ on $V^{U_x^{(e)}}$ is the same as
the $trace$ of $g\gamma' g^{-1}$ on $V^{U_{gx}^{(e)}}$.
\end{proof}

Let $g\in U^\circ$ and $\gamma'\in T^{\circ\circ}$. Since
$\gamma'$ and $g\gamma' g^{-1}$ have the same traces on fibers of
respective fixed points we see that the local constancy of the
expression ~(\ref{heartsuit}) holds for all such elements
$g\gamma' g^{-1}$. Since elements of this form contain an open
neighborhood of $\gamma$ we obtain the following analogue for
compact elements of Lemma ~\ref{ellopennbd}:
\begin{cor}\label{opennbd}
For $\gamma$ a (regular semisimple) compact element there is an
open neighborhood $\mathcal{O}_{U^\circ}(T^{\circ\circ})$ of
$\gamma$ in $G$ (which depends on $\gamma$ and $e$) such that for
all elements $\gamma'$ in this neighborhood the expression \emph{
~(\ref{heartsuit})} is constant.
\end{cor}

\subsection{Main result: the character on a compact element}

In this section we extend the Schneider-Stuhler formula to compact
elements. We start with a connected reductive group $G$ and a
finitely generated admissible representation $(\pi,V)$ of $G$,
with character function $\Theta_\pi(g)$. We assume that the
decomposition $(\oplus)$ is nice (independently of $e \geq e_0$).
We have shown the existence of such nice decompositions
(independent of $e\geq e_0$) for groups of semisimple rank $1$,
but in general we don't have the $e$ independence. Everything in
this section applies to a general $G$. It is only in the last
corollary, when we remove this assumption, that we will have to
restrict ourselves back to groups of semisimple rank $1$. Let
$\gamma \in G$ be a regular semisimple compact element, with $o
\in X^\gamma$, and fix $e \geq e_0$. The approach here is the
analogous approach taken in recovering the Schneider-Stuhler
result for semisimple rank $1$ groups for elliptic elements. [See
section ~\ref{SS:SL_2}].\\

Fix an open compact subgroup $K$ of $G$ with the following
properties (these properties are the analogue for compact elements
of the properties listed in Claim ~\ref{Kgroup}):
\begin{enumerate}
\item{} The character is locally constant on the
neighborhood $\gamma K$ of $\gamma$: $\Theta_\pi(\gamma)=
\Theta_\pi(\gamma K)$
\item{} $\gamma K$ is contained in the open neighborhood $U_{\gamma}^e$
of $\gamma$ described in Corollary ~\ref{opennbd}.
\item{} $\gamma$ normalizes $K$.
\end{enumerate}
As explained in section ~\ref{SS:SL_2}, it is possible to choose a
nice decomposition $(\oplus)$ which is common to all such
elements: the decomposition being $T_{\gamma k}$-equivariant
simultaneously for all $ \gamma \in K$ (by (3) above). Hence
similarly to the above we obtain:
$$trace(T_{\gamma k} Q_q^r;C_q) = trace(T_{\gamma
k} \overline{Q}_q^r;C_q) \quad d \geq q \geq 0 \qquad \mbox{for
all} \; k \in K. $$ Take $f:=\frac {1_{\gamma K}}{vol(K)} \in
C^{\infty}_c(G)$ to be the characteristic function of the set
$\gamma K$ normalized by its volume, and let $r_0=r_0(f)$ as in
Theorem ~\ref{thmfunction}. For all $r \geq r_0$ we obtain the
following:

\begin{eqnarray*}
\Theta_\pi(\gamma) & = & trace(\pi(f);V)\\
& = & \int\limits_K \sum_{q=0}^{d} (-1)^q trace(T_{\gamma
k}\overline{Q}_q^r; C_q)\,dk\\& = & \int\limits_K \sum_{q=0}^{d}
(-1)^q trace(T_{\gamma k}Q_q^r; C_q)\,dk\\ & = & \int\limits_K
\sum_{q=0}^{d} \underset{F(\gamma k)\in(X^{\gamma k} \cap
X^r)_q}\sum {(-1)}^q trace(\gamma k;V^{U_{F}^{(e)}})\,dk\\ & = &
\sum_{q=0}^{d} \underset{F(\gamma)\in(X^{\gamma} \cap X^r)_q}\sum
{(-1)}^q trace(\gamma;V^{U_{F}^{(e)}}).
\end{eqnarray*}

Which gives the main result of this paper:

\begin{thm}\label{main} Let $G$ be connected reductive, $(\pi, V)$ a finitely
generated admissible representation of $G$, $e \geq e_0(V)$, and
$\gamma$ regular semisimple compact. Assume the existence of nice
decompositions $(\oplus)$ uniformly in $e$ (see Definition
~\ref{def-nice-decomp}). For all $r$ large enough\footnote{ Let
$r'$ be large enough so that $X^{r'}$ contains a fundamental
domain $D$ for $X^\gamma$. Let $r_0=r_0(f)$ as above. The phrase
`$r$ large enough' means $r\geq max(r',r_0).$} we can express the
character $\Theta_\pi$ of $\pi$ on the compact element $\gamma$
using information contained in the truncated fixed point set
$X^\gamma \cap X^r$ as follows:
$$ \Theta_\pi(\gamma) = \sum_{q=0}^{d}
\underset{F(\gamma)\in(X^{\gamma} \cap X^r)_q}\sum {(-1)}^q
trace(\gamma;V^{U_{F}^{(e)}}) $$\\
\end{thm}

For groups of semisimple rank $1$, we know the existence of nice
decompositions uniformly in $e$ (Lemma ~\ref{key lemma}), so we
can remove this assumption in the theorem:
\begin{cor}\label{maincor}
For $G$ connected reductive of semisimple rank $1$, $(\pi, V)$ a
finitely generated admissible representation of $G$, $e \geq e_0$,
$\gamma$ regular semisimple compact, and $r$ large enough, the
following holds:
$$ \Theta_\pi(\gamma) = \sum_{q=0}^{d}
\underset{F\in(X^{\gamma} \cap X^r)_q}\sum {(-1)}^q
trace(\gamma;V^{U_{F}^{(e)}}) $$
\end{cor}

\section{Concluding remarks}
In Theorem ~\ref{main} we made the assumption about existence of
nice decompositions uniformly in $e$. Existence of nice
decompositions is controlled by the exactness of the truncated
complex $(\overline{*})$. In section ~\ref{nicedecomp} we proved
exactness of this complex: in the semisimple rank $1$ case, we
showed that exactness is independent of $e_r$, but in general we
had only shown exactness with dependence on $e_r$ (see Lemma
~\ref{key lemma}). Thus in the semisimple rank $1$ case, we are
able to remove the assumption and obtain Corollary ~\ref{maincor}.
One can hope that it might be possible to show exactness of
$(\overline{*})$ independently of $e_r$ and hence remove the
assumption from the main theorem.

\end{document}